\documentclass[leqno,11pt]{amsart} 
\usepackage{amsfonts,epsfig,latexsym,amssymb,amsmath,graphics,color,float,rotating,lscape,hyperref} 

\newtheorem{lemma}{Lemma}[section] 
 
\newtheorem{theorem}[lemma]{Theorem} 
\newtheorem{proposition}[lemma]{Proposition} 
\newtheorem{definition}[lemma]{Definition} 
\newtheorem{corollary}[lemma]{Corollary} 
\newtheorem{example}[lemma]{Example} 
\newtheorem{exercise}[lemma]{Exercise} 
\newtheorem{remark}[lemma]{Remark} 
\newtheorem{fig}[lemma]{Figure} 
\newtheorem{tab}[lemma]{Table} 
  
\newcommand{\bth}{\begin{theorem}} 
        \newcommand{\ethe}{\end{theorem}} 
\newcommand{\bre}{\begin{remark}\em } 
        \newcommand{\ere}{\end{remark}} 
\newcommand{\ble}{\begin{lemma}} 
        \newcommand{\ele}{\end{lemma}} 
\newcommand{\bde}{\begin{definition}} 
        \newcommand{\ede}{\end{definition}} 
\newcommand{\bco}{\begin{corollary}} 
        \newcommand{\eco}{\end{corollary}} 
  
\newcommand{\bpr}{\begin{proposition}} 
        \newcommand{\epr}{\end{proposition}} 
  
\newcommand{\bexer}{\begin{exercise}} 
        \newcommand{\eexer}{\end{exercise}} 
  
\newcommand{\bexam}{\begin{example}} 
        \newcommand{\eexam}{\end{example}} 
  
\newcommand{\bfi}{\begin{fig}} 
        \newcommand{\efi}{\end{fig}} 
  
\newcommand{\btab}{\begin{tab}} 
        \newcommand{\etab}{\end{tab}} 

\def\B_e{B_{\eta}(e)}

\definecolor{darkblue}{rgb}{.1, 0.1,.8} 
\definecolor{darkgreen}{rgb}{0,0.8,0.2} 
\definecolor{darkred}{rgb}{.8, .1,.1}

\newcommand{\beao}{\begin{eqnarray*}} 
        \newcommand{\eeao}{\end{eqnarray*}\noindent} 
  
\newcommand{\beam}{\begin{eqnarray}} 
\newcommand{\eeam}{\end{eqnarray}\noindent} 
  
\newcommand{\beqq}{\begin{equation}} 
\newcommand{\eeqq}{\end{equation}\noindent} 
  
\newcommand{\bce}{\begin{center}} 
        \newcommand{\ece}{\end{center}} 
  
\newcommand{\barr}{\begin{array}} 
        \newcommand{\earr}{\end{array}}

\newcommand{\vague}{\stackrel{\lower0.2ex\hbox{$\scriptscriptstyle 
                        \it{v} $}}{\rightarrow}} 
\newcommand{\dist}{\stackrel{\lower0.2ex\hbox{$\scriptscriptstyle 
                        \it{d} $}}{\rightarrow}} 
\newcommand{\weak}{\stackrel{\lower0.2ex\hbox{$\scriptscriptstyle 
                        \it{w} $}}{\rightarrow}} 
\newcommand{\what}{\stackrel{\lower0.2ex\hbox{$\scriptscriptstyle 
                        \it{\hat{w}} $}}{\rightarrow}} 
  
\newcommand{\bdis}{\begin{displaymath}} 
\newcommand{\edis}{\end{displaymath}\noindent}

\begin{document} 
        
  
       \title{Portfolio Selection under Multivariate Merton Model with Correlated Jump Risk}
        \author[B. Afhami]{Bahareh Afhami} 
        \author[M. Rezapour]{Mohsen Rezapour} 
        \author[M. Madadi]{Mohsen Madadi} 
                \author[V. Maroufy]{Vahed Maroufy} 
        \address[B. Afhami]{Department of Statistics, Faculty of Mathematics and Computer, Shahid Bahonar University of Kerman, Kerman, Iran.}\email{baharehafhami@math.uk.ac.ir} 
        \address[M. Madadi]{Department of Statistics, Faculty of Mathematics and Computer, Shahid Bahonar University of Kerman, Kerman, Iran.}\email{madadi@uk.ac.ir} 
        \address[M. Rezapour]{ 
                Department of Biostatistics \& Data Science, School of Public Health, University of Texas Health Science Center at Houston (UTHealth), 
                Houston, Texas, USA. }\email{Mohsen.R.Toughari@uth.tmc.edu} 
        \address[V. Maroufy]{ 
                Corresponding author: Department of Biostatistics \& Data Science, School of Public Health, University of Texas Health Science Center at Houston (UTHealth), 
                Houston, Texas, USA. }\email{Vahed.Maroufy@uth.tmc.edu} 
        
\begin{abstract} 
Portfolio selection in the periodic investment of securities modeled by a multivariate Merton model with dependent jumps is considered. The optimization framework is designed to maximize expected terminal wealth when portfolio risk is measured by  the Condition-Value-at-Risk ($CVaR$). Solving the portfolio optimization problem by Monte Carlo simulation often requires intensive and time-consuming  computation; hence a faster and more efficient portfolio optimization method based on closed-form comonotonic bounds for the risk measure $CVaR$ of the terminal wealth is proposed.
                
        \end{abstract} 
        
        \thanks{} 
        \maketitle 
  
\section{Introduction} \label{lll998}
The idea of portfolio optimization has been the focus of research groups and management teams since it was introduced by Markowitz in 1952 \cite{Markowitz:1952} and a wide variety of optimization methods using different financial models have been proposed. Dhaene et al. \cite{Dhaene:2005} studied multi-period portfolio selection  in a Black-Scholes market. Kaas et al. \cite{Kaas:2000} addressed the portfolio optimization problem  with a constant mix strategy using comonotonic approximations and Weert et al. \cite{Weert:2011} extended their results to the selection of a multi-period portfolio related to a specific provisioning problem. Brown et al. \cite{Brown:2011} considered dynamic portfolio optimization in a discrete-time finite-horizon setting, and Xu et al. \cite{Xu:2017} studied portfolio optimization for a Black-Scholes market with stochastic drift. Recently Wang et al. \cite{Wang:2016} implemented the Merton model to price a power exchange option. 

In this paper, we follow \cite{Wang:2016} and consider a multivariate Merton model with two kinds of dependent jumps. Uncommon jumps are such that the jump times in the marginal processes may differ from each other. These are individual jumps and occur upon the arrival of important new information  which, for example, may have a marginal effect on the price only of certain stocks or cause an imbalance between their supply and demand thus changing their prices.
Common jumps occur simultaneously in the marginal processes but the jump values may vary.
These simultaneous jumps can be caused by changes in the investment rate, alterations in the economic outlook or other reasons that impose a price jump at a specific time on all stock values.

\begin{figure}[ht] 
\centerline{\includegraphics[height=5 cm, width=6 cm]{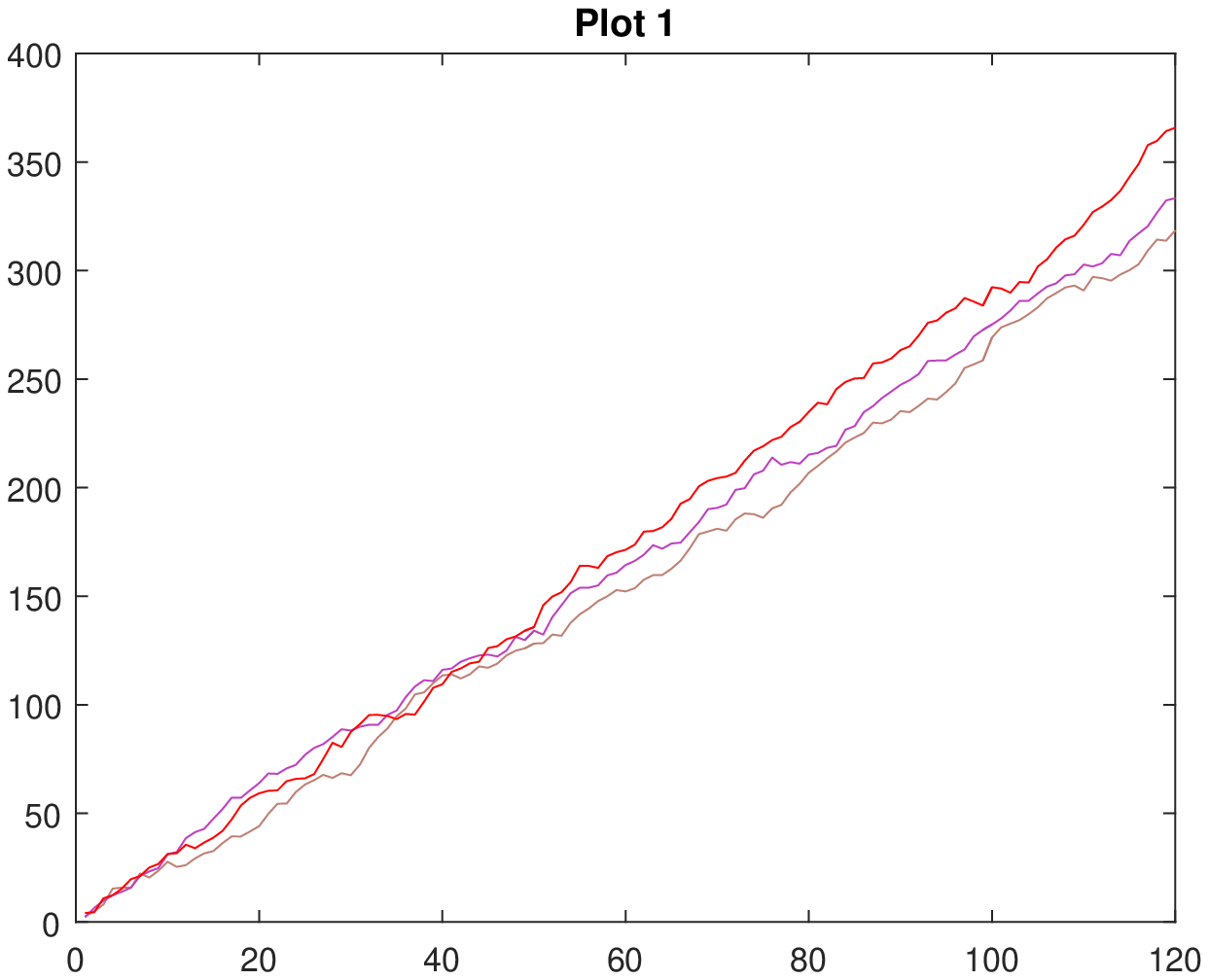} 
\includegraphics[height=5 cm, width=6 cm]{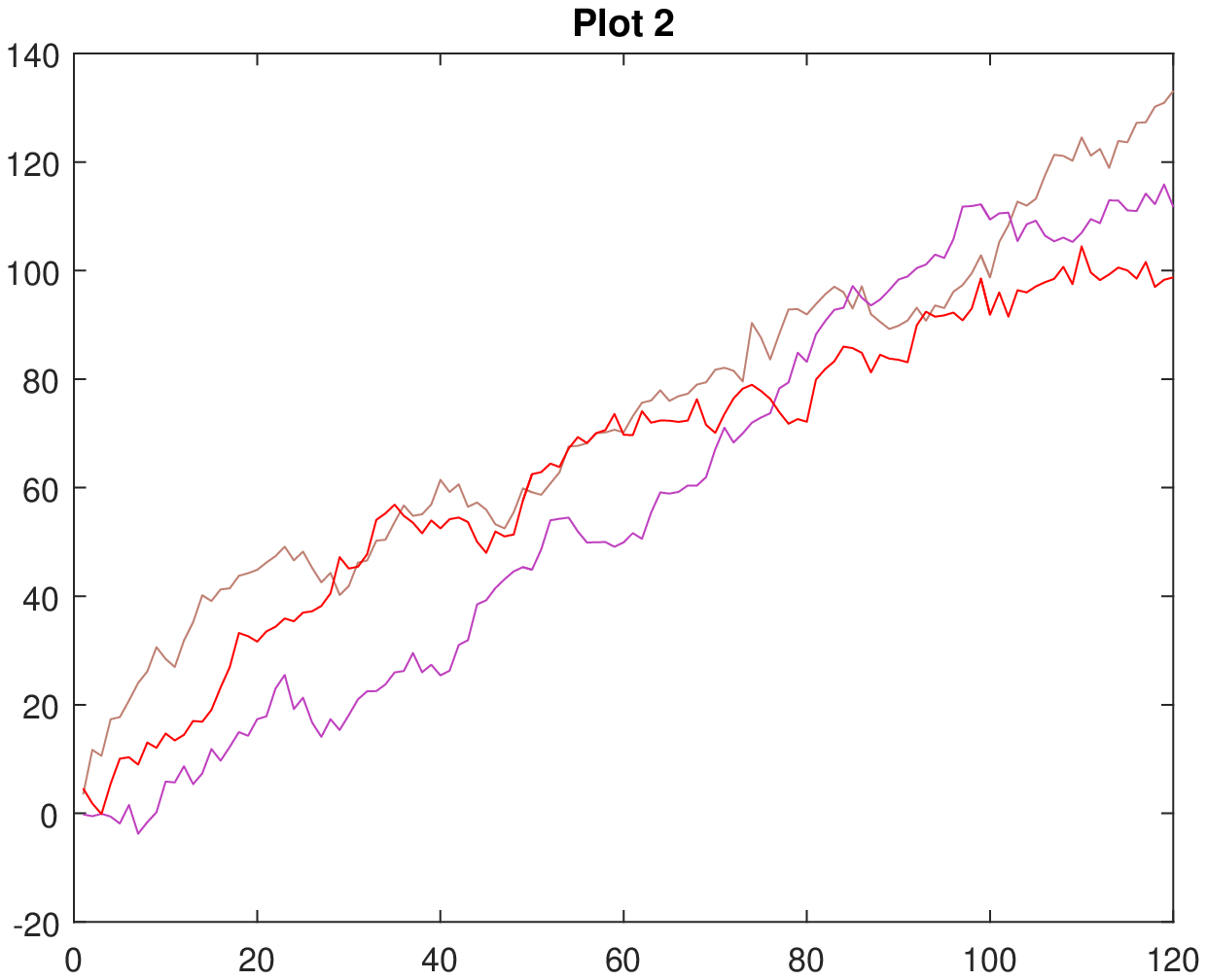}} 
\centerline{\includegraphics[height=5 cm, width=6 cm]{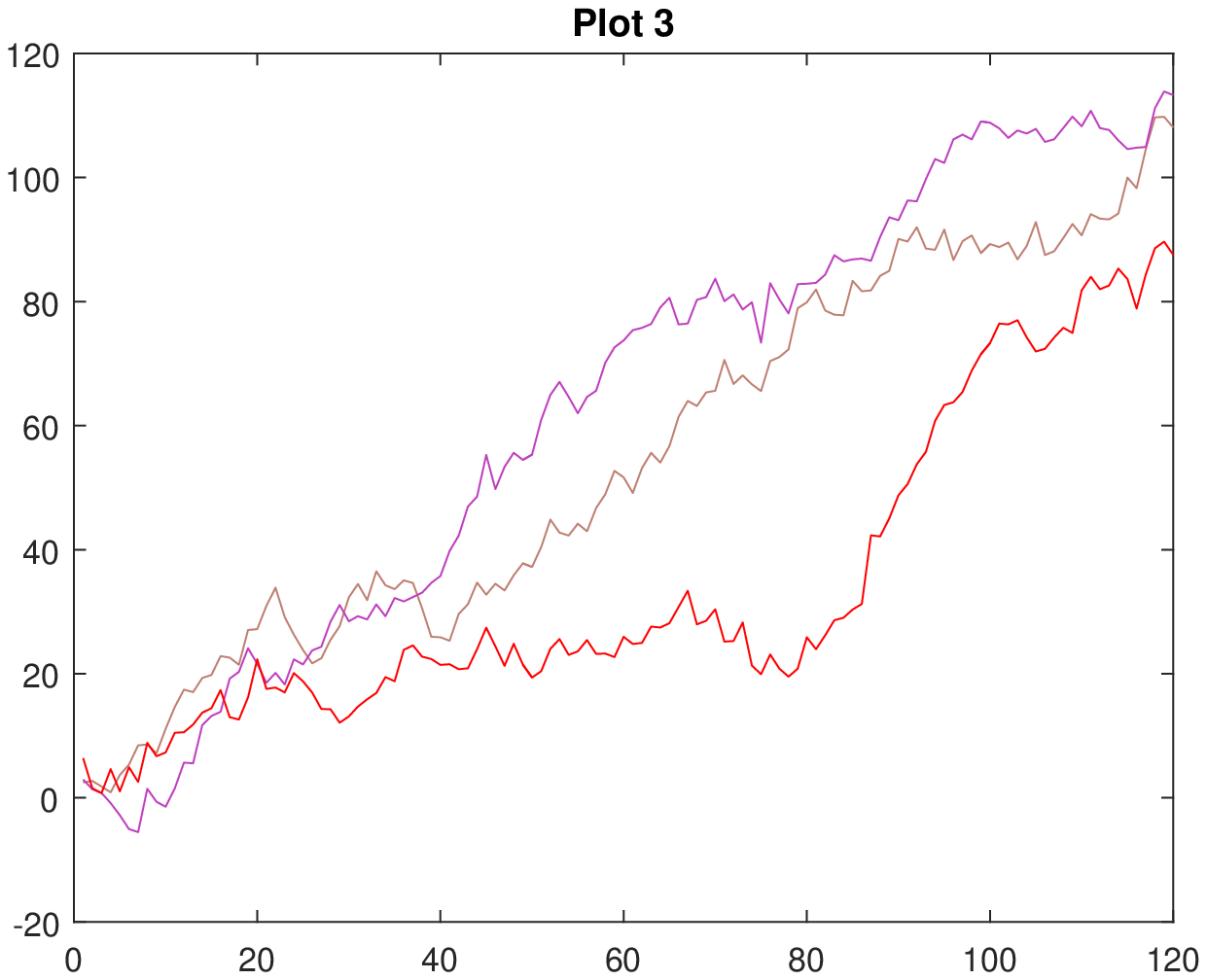} 
\includegraphics[height=5 cm, width=6 cm]{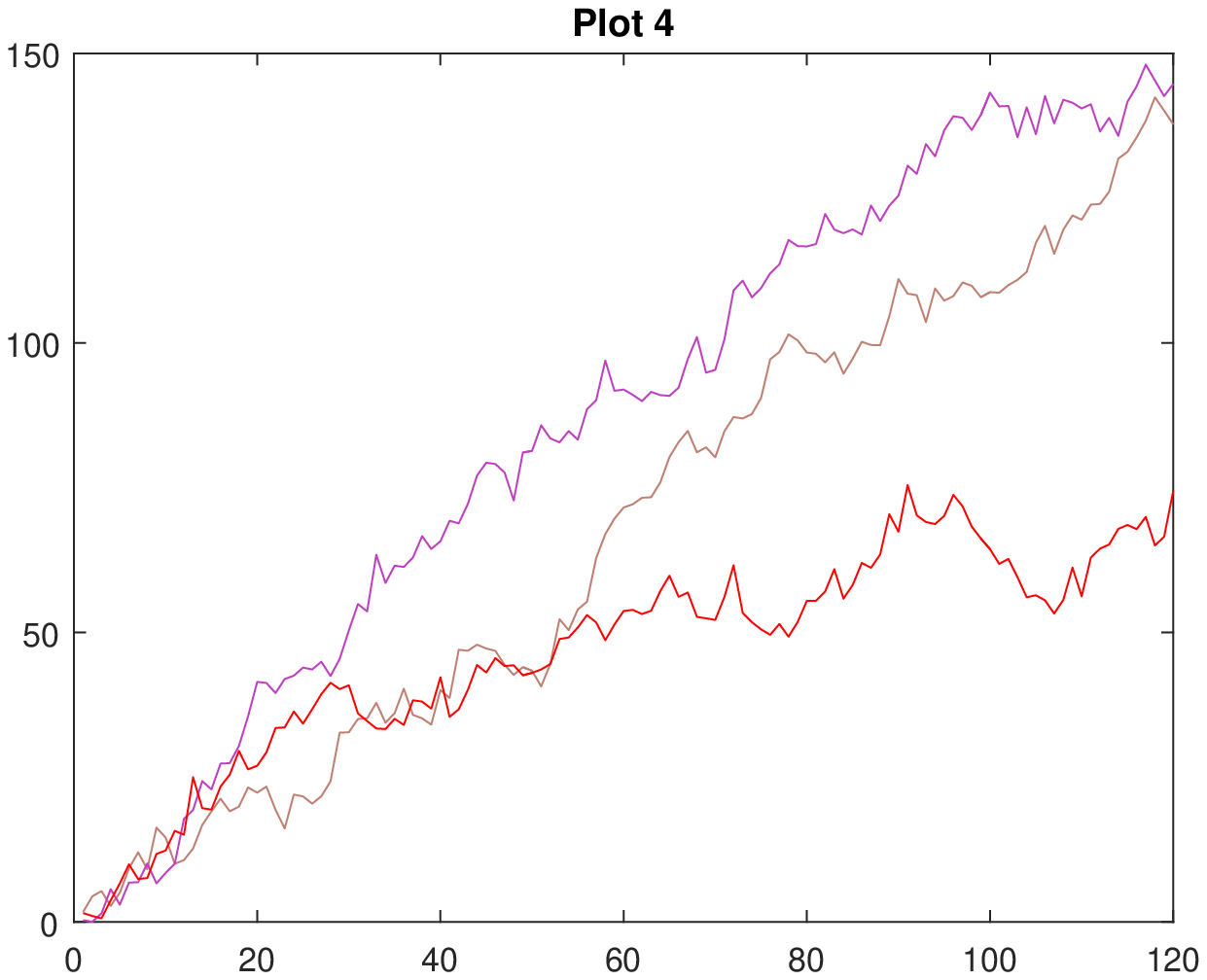}} 
\caption{\label{hgfffbgh} This plot presents four simulated processes: processes without any jump terms (Plot 1), processes with common jumps but maybe different jump values (Plot 2), processes without common jumps (Plot 3), and processes with both types of jumps (Plot 4).} 
\end{figure}

The flexibility of the Merton model enables it to fit financial data with correlated jumps and can be used to optimize a portfolio for a periodic investment problem in a constant mix strategy framework. Figure \ref{hgfffbgh} illustrates this flexibility and compares it with three other models using simulated processes. The figure presents simulated data from four models, including processes without any jump terms (Plot 1), processes with common jumps but with different jump values (Plot 2), processes without common jumps (Plot 3), and processes with both types of jumps (Plot 4). These plots show the ability of the Merton model to take into account both individual jumps as well as jumps occurring at a specific times. 

In order to solve the optimization problem for this model we need to solve a stochastic differential equation (SDE) corresponding to a portfolio containing multiple assets. In what follows we provide a closed-form solution for the terminal wealth process of the SDE corresponding to a constant mix strategy investment problem. To optimize a portfolio we maximize the terminal wealth process by controlling its risk measure $CVaR$ and rather than applying numerical methods such as Monte Carlo Markov Chain (MCMC) to solve this optimization problem, we propose an alternative method using {\it comonotonic} bounds. This approach reduces computational complexity and returns a closed-form solution to the portfolio optimization problem for maximizing the expectation of terminal wealth when the corresponding risk measure of the comonotonic bound is under control. 

The concept of comonotonicity has been thoroughly studied in actuarial and financial literature (Dhaene et al. \cite{Dhaene:2002a} and Dhaene et al. \cite{Dhaene:b}) and the assumption of comonotonic dependence structures between asset prices has proved to be a useful when solving financial problems such as derivative pricing (Deelstra et al. \cite{{Deelstra:2004},{Deelstra:2008}} and Linders et al. \cite{Deelstra:2016}). Vanduffel et al. \cite{Vanduffel:2008} used comonotonic approximations to obtain closed forms for the lower and upper bounds for the price of a continuously sampled European-style Asian option with fixed exercise price.

The rest of the paper is organized as follows. In Section \ref{lll998} we introduce the Merton model,  and review existing definitions and preliminary results. In Section \ref{Main results}, we present the main results; we first solve the SDE corresponding to the terminal wealth process of a portfolio following a jump-diffusion model containing one risk-free asset and $m$ risky assets and then obtain comonotonic bounds for the terminal wealth in order to solve the portfolio optimization problem while controlling $CLVaR$ of the terminal wealth. Proofs of all the results stated in this section may be found in the appendix. In Section \ref{Numerical illustration}, we carry out our optimization method on simulated data  and data containing the daily prices of the stocks of Zoom Video Communication Inc. (ZM), and Tesla Inc. (TSLA) obtained from the $S\&P\, 500$ market index. 
  
\section*{Preliminaries} 
The base asset price in a Black-Scholes model follows the geometric Brownian motion process and the dynamics of the stock price are described by a continuous-time diffusion process whose sample path is continuous with probability one. Merton \cite{Merton1973} illustrated that in order to satisfy the conditions of the Black-Scholes model, trading should take place continuously in time and the price dynamics of the stock should have a continuous sample path with probability one. 
However, since real world datasets are discrete realizations of variables and a continues path model does not provide for jumps, it is generally not suitable for real stock prices. 

Following \cite{Wang:2016}, we shall assume that all continuous time processes are Wiener and jump components are Poisson driven. 
Merton's jump-diffusion model  \cite{Merton1} is commonly used in literature and has abundant applications in financial data analysis. For example, Kolmanovsky et al. \cite{Kolmanovsky:2002} applied this model to address the optimal containment control problem in a financial market,  Wang et al. \cite{Wang:2017} investigated a power exchange option pricing problem assuming a jump-diffusion model, and  Xu et al. \cite{Xu:2019}  provided an  analytical valuation  of power exchange options based on the Merton model with dependent jumps. The model can be written as 
$$\frac{dP}{P}=(\alpha-\lambda\, k)dt+\sigma\, dB(t)+ZdN(t), $$ 
where $\alpha$ is the instantaneous expected return on the stock, $\sigma^2$ is the instantaneous variance of the return conditional on arrival of no important new information, $dB$ is a stochastic differential of the standard Brownian motion, and $dN$ is a stochastic differential of a Poisson process defined based on the natural filtration of the random processes $B$ and $N$. Moreover, $dN$ and $dB$ are assumed to be independent. 
  
Bayraktar et al. \cite{Bayraktar:2008} illustrated the advantages of the Merton model over a piecewise Markov model and Giesecke et al. \cite{Giesecke:2018} used it to study asset, commodity and energy prices, interest and exchange rates, and the timing of corporate and sovereign defaults. 
  
In this paper, we study the return of a portfolio consisting of one risk-free asset together with several risky assets which follow the jump-diffusion Merton model. To be more specific, we assume that the jump risk common to all assets is  diversifiable in the market and has a zero risk premium. This assumption ensures that there exists a risk neutral measure under which (for $j=1,2,\ldots,m$), the dynamics of risky assets prices become 
\begin{eqnarray} \label{equation1} \nonumber 
\frac{dP_j(t)}{P_j(t-)}&=&(r+\mu_j -\lambda\, h_{j,0}-\lambda_j\, h_{j,1})\, dt+\sigma_j\, dB'_j(t)+(e^{Z_{\iota(t-),j,0}}-1)\, dN(t)\\ 
&&+(e^{Z_{\iota(t-),j,1}}-1)\, dN_j(t), 
\end{eqnarray} 
where $\mu_j$ is the drift of geometric Brownian motion, $\sigma_j>0$ is the volatility of the asset $P_j\,( j=1,2,\ldots,m)$ and $B'_j$ are standard Brownian motions with  $Cov(B'_i(t),B'_j(t+s))=t\,\rho_{i,j}\, (t,s\geq 0,\,   i,j=1,\ldots,m)$ and $h_{j,\ell}=E[e^{Z_{\iota(t),j,\ell}}]-1\, (\ell=0,1).$ Moreover, $N(t)$ and $N_j(t)$ are independent Poisson processes with intensity rates $\lambda$ and $\lambda_j$, respectively. Note that the SDE \eqref{equation1} is defined on the natural filtration of the random processes $(B'_1,\ldots,B'_m)$ and $N,N_1,\ldots,N_m\,.$ We will use  Theorem 6 (pg. 249) in \cite{Protter} in the proof of Proposition \ref{final1} to show there exits a probability space on which the solution to the SDE \eqref{equation1} is unique. Discounted changes in assets price are divided into two groups: individual changes corresponding to $N_j(t)$ and common changes corresponding to $N(t).$ As in  \cite{Wang:2016}, for $j=1,2,\ldots,m$ and $k=1,2,\ldots,N(t),\, Z_{k,j,0}$ denotes the jump magnitude of the $k$-th common jump for the asset $P_j$ in $(0,t]$ while for $k=1,2,\ldots,N_j(t),\, Z_{k,j,1}$ denotes the $k$-th individual jump of asset $P_j,$ in $(0,t].$  We also assume that for all $j,\, k,$ and $\ell=0,1$ the random variables $Z_{k,j,\ell}$ are independent of $B'_j(t)$ and for fixed $j=1,2,\ldots,m,$ the random variables $Z_{k,j,\ell}$ are $i.i.d$ for all $k$ and $\ell.$ The variables $Z_{k,j,\ell}$ and $Z_{k,j',\ell}$ are also independent for $j\neq j',$ but they may not be identically distributed. 
  
We close this section by recalling definitions and theorems which will be required in Section \ref{Main results}. 
\bde {\bf(Comonotonicity)} 
The random vector $(X_1, X_2, \ldots, X_n)$ is comonotonic if 
\begin{eqnarray*} 
 (X_1, X_2, \ldots, X_n)\sim (F_1^{-1}(U), F_2^{-1}(U), \ldots, F_n^{-1}(U)), 
\end{eqnarray*} 
where $\sim$ means equality in distribution, $F_k^{-1}$ is the (generalized) inverse distribution function of $X_k$ for $k = 1, 2, \ldots, n,$ and $U$ is a uniform random variable on $(0,1).$ 
\ede 
\bde {\bf(Convex order)}. 
A random variable $X$ is said to precede $Y$ in the sense of convex order $(X \leq_{cx} Y)$, if $\mathbb{E}(X) =\mathbb{E}(Y)$ and $\mathbb{E}[\max(X-d , 0)] \leq \mathbb{E}[\max(Y-d , 0)]$ for all real $ d.$ 
\ede 
The following fundamental theorem (Kaas et al. \cite{Kaas:2000}) provides convex bounds for sums of random variables. 
\bth\label{theorem1.4} 
{\bf(Convex bounds for sums of random variables)}. For any random vector 
$(X_1,X_2,\ldots,X_n)$ and any random variable $\Lambda$, 
\begin{eqnarray*} 
\sum\limits_{j=1}^{n}\mathbb{E}(X_{j}\mid \Lambda) \leq_{cx} \sum\limits_{j=1}^{n}X_{j} \leq_{cx} \sum\limits_{j=1}^{n}F_{X_{j}}^{-1}(U). 
\end{eqnarray*} 
\ethe 

In this paper we consider the risk measures $VaR,\, CVaR$ and $CLVaR.$ Calculating these risk measures for sums of random variables is not simple, but accurate and easy comonotonic approximations are available. 

\bde {\bf($VaR, CVaR, CLVaR$)} 
For the set of real numbers $\mathbb{R}$ and $p\in (0,1),$\\ 
(i) The Value-at-Risk 
\begin{eqnarray*} 
VaR_{p}(X):=F_{X}^{-1}(p)=\inf \{ x\in \mathbb{R} \mid F_{X}(x) \geq p \}, 
\end{eqnarray*} 
where $F_{X}(x)=Pr(X \leq x) $ and by convention, $\inf \{\emptyset \} = + \infty.$\\ 
(ii) The Conditional-Value-at-Risk
\begin{eqnarray*} 
CVaR_{p}(X):=\mathbb{E} \left( X \mid X >VaR_{p}(X) \right). 
\end{eqnarray*} 
(iii) The Conditional-Left-side-Value-at-Risk 
\begin{eqnarray*} 
\hspace{.2cm} CLVaR_{p}(X):=\mathbb{E} \left( X \mid X < VaR_{p}(X) \right). 
\end{eqnarray*} 
\ede

It is obvious that 
\begin{eqnarray}\label{relation1} 
\hspace{-.8cm} CVaR_{1-p}(X)=-CLVaR_{p}(-X). 
\end{eqnarray} 

\bth\label{Proposition1} 
(Dhaene et al. \cite{Dhaene:2005}) If the random vector $(X_1,X_2,\ldots,X_n)$ is comonotonic and $S=X_1+X_2+\ldots+X_n$, for all $p\in(0,1),$ then
\begin{eqnarray*} 
VaR_p(S)=\sum\limits_{j=1}^n VaR_p(X_j). 
\end{eqnarray*} 
Also, if all marginal distributions $F_k$ are continuous, then
\begin{eqnarray*} 
CVaR_p(S)=\sum\limits_{j=1}^nCVaR_p(X_j). 
\end{eqnarray*} 
\ethe 
  
\section{Main results}\label{Main results} 
In this section, we consider the terminal wealth of a portfolio modeled by the Merton model  (\ref{equation1}) and obtain a closed form expression for the terminal wealth of periodic investment on a portfolio with risk-free as well as risky assets. For this purpose, we consider a market where $(m + 1)$ assets, (one risk-free and $m$ risky) are traded continuously. Suppose that at time $t,$ the decision maker invests a fraction $x_j(t) \,(j=1,2,\ldots,m)$ of his/her wealth in the $j$-th risky asset  and $1-\sum\limits_{j=1}^m x_j(t)$ in the risk-free asset. Now let $\alpha_0, \alpha_1,\ldots,\alpha_{\tau-1}$ denote periodic endowments at predetermined points $0,1,\ldots,\tau-1$ and suppose that the decision maker can only rebalance his/her portfolio at the beginning of each period.   Also suppose that $P_j(t),$ the price of the $j$-th risky asset at time $t,$ satisfies (\ref{equation1}) and $Z_{k,j,\ell}\, (j=1,\ldots,m,\, \ell=0,1)$ are continuous random variables with $\mathbb{E}(Z_{k,j,\ell})=\mu_{Z_{1,j,\ell}}$ and $Var(Z_{k,j,\ell})=\sigma^2_{Z_{1,j,\ell}}.$ 

Assume that $P_0(t),$ the price of the risk-free asset at time $t,$ satisfies 
\begin{eqnarray}\label{equation1jgvkj} 
\frac{dP_0(t)}{P_0(t- )}=r\, dt, \hspace{0.5cm} r>0. 
\end{eqnarray} 
$ {\mathbf{\mu}}:=(r+\mu_1-\lambda h_{1,0}-\lambda_1 h_{1,1},\ldots,r+\mu_m-\lambda h_{m,0}-\lambda_m h_{m,1})^\top,$   ${\mathbf {x}}:=(x_1,\ldots,x_m)^\top,$  
$\mu({\mathbf{x}}):=({\mathbf{\mu}}-r {\mathbf{1}})^\top\, {\mathbf{x}}+r,\,$ 
$B_{{\mathbf{x}}}(t):=\frac{1}{{\mathbf{x}}^\top\,{\mathbf{\Sigma}}\,\, {\mathbf{x}}}\sum\limits_{j=1}^{m} x_j\, \sigma_j\, B'_j(t),$ and ${\mathbf{\Sigma}}_{m\times m}:={\left(\sigma_{i}\,\sigma_{j}\,\rho_{i,j}\right)}_{i,j}.$ Then, according to (\ref{equation1}) and (\ref{equation1jgvkj}), the return  within the $k$-th period is given by
\begin{eqnarray}\label{stochastic equation}\nonumber 
\frac{dP(t)}{P(t-)}&=&\sum\limits_{j=1}^{m}x_j\, \frac{dP_j(t)}{P_j(t-)}+(1-\sum\limits_{j=1}^{m}x_j) \frac{dP_0(t)}{P_0(t-)}\\ \nonumber 
&=&\mu\left({\mathbf{x}}\right)\,dt+\left({\mathbf{x}}^\top\,{\mathbf{\Sigma}}\,\, {\mathbf{x}}\right)\,dB_{{\mathbf{x}}}(t)\\ 
&&+\sum\limits_{j=1}^{m} x_j\, (e^{Z_{\iota(t-),j,0}}-1)\,dN(t)+\sum\limits_{j=1}^{m} x_j\,(e^{Z_{\iota(t-),j,1}}-1)\,dN_j(t). 
\end{eqnarray} 
The following proposition provides a solution to the SDE given in (\ref{stochastic equation}). 
  
\bpr\label{final1} 
Suppose the random processes $N,N_{j_1}$ and $N_{j_2}$ do not have common jumps for $j_1\neq j_2=1,2,\ldots,m.$  Then the solution of (\ref{stochastic equation}) is 
\begin{eqnarray*} 
P(t)&=& e^{\left(\mu({\mathbf{x}})-\frac{{\mathbf{x}}^\top {\mathbf{\Sigma}}\, {\mathbf{x}}}{2}\right)\,t+({\mathbf{x}}^\top{\mathbf{\Sigma}}\, {\mathbf{x}})\, B_{{\mathbf{x}}}(t)+\sum\limits_{j=1}^{m}  \sum\limits_{k=1}^{N(t)}Z^\star_{k,j,0}({\mathbf{x}})+\sum\limits_{j=1}^{m} \sum\limits_{k=1}^{N_j(t)}Z^\star_{k,j,1}({\mathbf{x}})}, 
\end{eqnarray*} 
where $e^{Z^\star_{k,j,\ell}({\mathbf{x}})}-1=x_j(e^{Z_{k,j,\ell}}-1),\, \ell=0,1,\, j=1,\ldots,m.$ 
\epr 
  
  
We note that an investment of a unit amount of wealth at time $t-1$ will grow to $e^{Y_{t}({\mathbf{x}})}$ at time $t,$ where{\small 
\begin{eqnarray*}\label{ypil}\nonumber 
Y_{t}({\mathbf{x}})&=&\mu({\mathbf{x}})-\frac{{\mathbf{x}}^\top{\mathbf{\Sigma}}\,{\mathbf{x}}}{2}+\sum\limits_{j=1}^{m}  \sum\limits_{k=N(t-1)+1}^{N(t)}Z^\star_{k,j,0}({\mathbf{x}}) 
+\sum\limits_{j=1}^{m} \sum\limits_{k=N_j(t-1)+1}^{N_j(t)}Z^\star_{k,j,1}({\mathbf{x}})\\ 
&&+({\mathbf{x}}^\top{\mathbf{\Sigma}}\, {\mathbf{x}})\left(B_{{\mathbf{x}}}(t)-B_{{\mathbf{x}}}(t-1)\right). 
\end{eqnarray*}} 
So the wealth $W_t$ at the end of the $t$-th period will satisfy the recursion equation 
\begin{eqnarray*}  
W_t=W_{t-1}\,e^{Y_{t}({\mathbf{x}})}+\alpha_t, \hspace{.75cm} t=1,2,\ldots, \tau, 
\end{eqnarray*} 
with $W_0=\alpha_0.$ Hence the terminal wealth can be written as 
\begin{eqnarray}\label{kjjjjjjj} 
\hspace{-3.5cm} W_\tau=\sum\limits_{t=0}^{\tau-1}\alpha_t\, e^{S_t+V_t}, 
\end{eqnarray} 
where 
\beao 
&&S_t=(\tau-t)(\mu({\mathbf{x}})-\frac{{\mathbf{x}}^\top{\mathbf{\Sigma}}\, {\mathbf{x}}}{2}) 
+\sum\limits_{j=1}^{m}  \sum\limits_{k=N(t)+1}^{N(\tau)}Z^\star_{k,j,0}({\mathbf{x}}) 
+\sum\limits_{j=1}^{m} \sum\limits_{k=N_j(t)+1}^{N_j(\tau)}Z^\star_{k,j,1}({\mathbf{x}})\,,\\ 
&&V_t=({\mathbf{x}}^\top{\mathbf{\Sigma}}\, {\mathbf{x}}) \left(B_{{\mathbf{x}}}(\tau)-B_{{\mathbf{x}}}(t)\right). 
\eeao 
  
Due to limitation in borrowing from the risk-free asset, the proportion ${\mathbf{x}}$ has an upper bound and so $\mu({\mathbf{x}})$ has an upper bound. We denote this upper bound by $c_0.$ Therefore the terminal wealth problem reduces to the following optimization problem: 
\begin{eqnarray}\label{aim1} 
&&\max\limits_{\mathbf{x}}\, \mathbb{E}(W_\tau),\\ \nonumber 
&&CLVaR_{p}(W_\tau)\geq K,\\ \nonumber 
&&\mu({\mathbf{x}}) \leq c_0. 
\end{eqnarray} 
Although, the optimization problem (\ref{aim1}) can be solved using MCMC approach, it is inefficient and slow, especially when $\tau$ and $m$ are large. Also, in general, the exact distribution of $W_\tau$ cannot be obtained. Therefore we search for the solution to this optimization problem using comonotonic bounds. We first proceed to calculate the comonotonic lower bound. 
\bth\label{Gnllllnnn} 
The comonotonic lower bound for terminal wealth for fixed $\Lambda=\sum\limits_{t=0}^{\tau-1} \alpha_t\, V_t $ is 
\begin{eqnarray}\label{termG_n^lnnnnn} 
\hspace{-2.7cm} W_\tau^L=\sum\limits_{t=0}^{\tau-1}\alpha_t\, e^{\frac{\Lambda}{\sigma_\Lambda}\, \sqrt{{\mathbf{x}}^\top{\mathbf{\Sigma}}\, {\mathbf{x}}}\,\, c_{4,t}+c_{3,t}({\mathbf{x}})+c_{2,t}}, 
\end{eqnarray} 
where $c_{2,t},\, c_{3,t}$ and $c_{4,t}$ are given (\ref{a14}), (\ref{a15}) and (\ref{a16}).  
\ethe 
  
It is common to use the first-order Taylor expansion of the exponential function in (\ref{termG_n^lnnnnn}) to obtain a closed-form solution to portfolio optimization (see for example \cite{Xu:2017}).   
\bco\label{corollary3.77} 
An approximation of the comonotonic lower bound for terminal wealth is as follows: 
\begin{eqnarray}\label{h34567} 
W_\tau^{'L}=c_{5}+ c_{6}\, ({\mathbf{\mu}}-r {\mathbf{1}})^\top\, {\mathbf{x}}-c_7\, ({\mathbf{x}}^\top {\mathbf{\Sigma}}\, {\mathbf{x}})+c_8\, \frac{\Lambda}{\sigma_\Lambda}\, \sqrt{{\mathbf{x}}^\top {\mathbf{\Sigma}}\, {\mathbf{x}}}, 
\end{eqnarray} 
where $c_{5},\, c_{6},\, c_{7}$ and $c_{8}$ are given in (\ref{a17}), (\ref{a18}), (\ref{a19}) and (\ref{a20}). 
\eco 
  
In order to solve problem (\ref{aim1}) we first plug in $W_\tau^{'L}$ for $W_\tau.$ In other words, we first solve the following problem  
\begin{eqnarray}\label{n88} 
&&\max\limits_{\mathbf{x}}\, \mathbb{E}(W_\tau^{'L}),\\ \nonumber 
&&CLVaR_{p}(W_\tau^{'L})\geq K,\\ \nonumber 
&&\mu({\mathbf{x}}) \leq c_0.\nonumber 
\end{eqnarray} 
From (\ref{relation1}) it follows that the condition $CLVaR_{p}(W_\tau^{'L})\geq K$ is equivalent to 
$CVaR_{1-p}(-W_\tau^{'L})\leq -K.$

Now it can be easily shown that 
\begin{eqnarray*}\label{expwtl2222} 
\hspace{-4cm}\mathbb{E}(W_\tau^{'L})=c_{5}+ c_{6}\, ({\mathbf{\mu}}-r {\mathbf{1}})^\top\, {\mathbf{x}}-c_7\, ({\mathbf{x}}^\top {\mathbf{\Sigma}}\, {\mathbf{x}}).
\end{eqnarray*} 
We may also compute also the expression for $CLVaR_{p}(W_\tau^{'L})$ according to the  following lemma. 
\ble \label{gghhjk54} 
The Conditional-Value-at-Risk of $W_\tau^{'L}$ equals 
\begin{eqnarray}\label{cvarwtl} 
CVaR_{1-p}(-W_\tau^{'L})=-c_5-c_6\, ({\mathbf{\mu}}-r {\mathbf{1}})^\top\, {\mathbf{x}}+c_7\, ({\mathbf{x}}^\top {\mathbf{\Sigma}}\, {\mathbf{x}}) 
+c_{9}\,\sqrt{{\mathbf{x}}^\top {\mathbf{\Sigma}}\, {\mathbf{x}}}, 
\end{eqnarray} 
where $c_{9}$ is given in (\ref{a21}).  
\ele 

A constant mix strategy implies keeping the investment proportions constant. In other words, the investor may buy and sell stocks along with changing market conditions, but the proportions of the various assets in the portfolio must remain constant. In the rest of this section we solve problem (\ref{n88}) for the constant mix strategy by using fractional the Kelly-strategy considered by MacLean et al. \cite{MacLean:2006}. 
  
\bpr\label{final4} The solution of problem (\ref{n88}) for $p<0.5$ is of the form ${\mathbf{x}}=q\,{\mathbf{x}}^\star$ where ${\mathbf{x}}^\star={\mathbf{\Sigma}}^{-1}\, ({\mathbf{\mu}}-r {\mathbf{1}}),$ and $q$ is as in the next proposition. 
\epr 
The following results can be shown using  similar arguments to the proof of Proposition 3 in  Xu et al. \cite{Xu:2017}.
  
\bpr\label{final6} 
The value of $q$ in Proposition \ref{final4} is of the form $q=\min\{q_1,\, q_2,\, q_3\},$ with 
$q_1=\frac{-B_2-\sqrt{B_2^2-4\, B_1\, B_3} }{2\, B_1},\, q_2=\frac{c_6}{2\, c_7},$ and 
$ q_3=\frac{c_0-r}{({\mathbf{\mu}}-r {\mathbf{1}})^\top {\mathbf{\Sigma}}^{-1}\,({\mathbf{\mu}}-r {\mathbf{1}})},$ where 
\begin{eqnarray*} 
B_1=-c_7\,\, ({\mathbf{\mu}}-r {\mathbf{1}})^\top\, {\mathbf{\Sigma}}^{-1}\,({\mathbf{\mu}}-r {\mathbf{1}}),\hspace{.7cm} 
B_2=\frac{c_6}{-c_7}\, B_1-c_9\, \sqrt{\frac{B_1}{-c_7}},\hspace{.7cm} 
B_3=c_5-K, 
\end{eqnarray*} 
and $c_{5},\, c_{6},\, c_{7}$ and $c_{9}$ are given in (\ref{a17}), (\ref{a18}), (\ref{a19}) and (\ref{a21}). 
\epr 

\section{Simulation study and Data Example }\label{Numerical illustration} 
In Example \ref{example3.11} we compute the optimal constant mix strategy for simulated data. In Example \ref{example3.33} besides computing the optimal constant mix strategy for data extracted from the $S\&P\, 500$ stock market index, we also compute the optimal value of terminal wealth.  

\begin{figure}[ht] 
\centerline{\includegraphics[width=9cm]{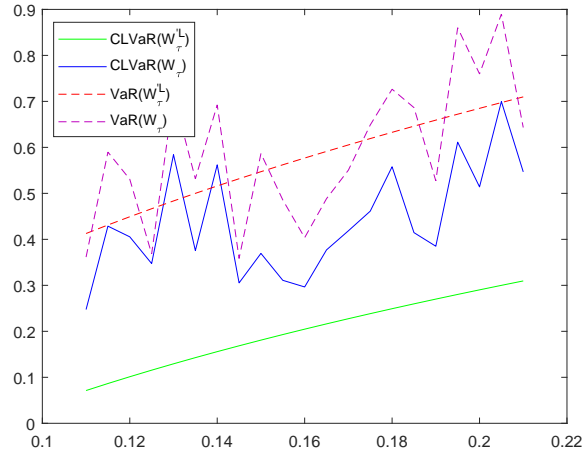}} 
\caption{\label{jjj8876} Diagram of $VaR$ and $CLVaR$ and their corresponding comonotonic approximations for different proportion of risk budget. } 
\end{figure}

 To assess the accuracy of comonotonic approximations for estimating $VaR$ and $CVaR$, we consider the risk measures of a random process satisfying \eqref{equation1}, and for simplicity we simulate a univariate random process (i.e., $m=1$) with $ \tau=1,\, r=0.03,\, \mu_j=\lambda=\lambda_1=\sigma_1^2=0.5,\,\mu_{Z_{1,1,\ell}}=\sigma_{Z_{1,1,\ell}}^2=0.1\, (\ell=0,1) $ and calculate the values of $VaR(W),\, CLVaR(W),\, VaR(W_\tau^{'L})$ and $CLVaR(W_\tau^{'L}),$  for different values of $p.$
Figure \ref{jjj8876} displays the original values of the risk measures $VaR$ and $CVaR$  and their corresponding  estimated values. 

\bexam \label{example3.11} 
Suppose $Z^\star_{k,j,\ell}({\mathbf{x}})$ for $k \in {\mathbb{Z}_{+}},\, j=1,\ldots,3,\,$ and $\ell=0,1,$ have a normal distribution with mean and variance $\mu_{Z_{1,j,\ell}},\sigma_{Z_{1,j,\ell}}^2,$ respectively. We simulate data from a normal distribution with the following parameters: 

\begin{eqnarray*} 
&&\tau=3,\,\,\, m=3,\,\,\, r=0.03,\,\,\, p=0.05,\,\,\,  \alpha_t=(1,1,1,1)^\top,\, \\ 
&& \lambda=0.9,\,\,\, \lambda_j=(0.82,0.8,0.81)^\top,\,\,\, \mu_j=(1.5,1.5,1.5)^\top \\ 
&&\sigma^2_{Z_{1,j,\ell}}=(0.063,0.062,0.061)^\top,\,\,\, \mu_{Z_{1,j,\ell}}=(0.041,0.042,0.043)^\top,\, \ell=0,1, 
\end{eqnarray*} 
and 
\begin{displaymath} 
{{\mathbf{\Sigma}}}=\left( \begin{array}{ccc} 
{1.3689} & {1.3455} & {1.3501}  \\ 
{1.3455} & {1.3689} & {1.3501}  \\ 
{1.3501} & {1.3501} & {1.3877} \\ 
\end{array} \right). 
\end{displaymath}  
In Equation (\ref{n88}) we set 
$K=k^\star\, \sum\limits_{i=1}^{\tau} \alpha_j e^{\frac{(\tau-i+1)r}{\tau}}$
 where $k^\star$ is the stop loss rate. The decision maker using a constant mix strategy invests a fraction $x_j$ of his/her capital in the risky asset $P_j.$ The optimal values of the portfolio for different values of $k^\star$ are given in Table \ref{mmoi90}. This table illustrates that  $\sum\limits_{j=1}^{m}x_j $ decreases as $k^\star$ increases. 
 \begin{table}[h!] 
 \center \caption{Constant mix strategy.} \label{mmoi90} 
 \begin{tabular}{ c c c  } 
    \noalign{\smallskip} \hline 
    $k^\star$ & $ (x_1,x_2,x_3)$ & $Total$   \\  \hline  
    $0.5$ & $(0.3375,0.3622,0.2920)$ & $0.9916$   \\ 
    $0.6$ & $(0.3233,0.3469,0.2797)$ & $0.9498$  \\  
    $0.7$ & $(0.3085,0.3310,0.2669)$ & $0.9063$   \\ 
    $0.8$ & $(0.2930,0.3144,0.2535)$ & $0.8609$   \\ 
    $0.9$ & $(0.2768,02970,0.2394)$ & $0.8132$   \\ \hline                               
 \end{tabular} 
 \end{table} 
  
\eexam 
  
\bexam \label{example3.33} 

In this example we consider stocks of Zoom Video Communication Inc. stock (ZM), an American communications technology company with headquarters in San Jose, California and Tesla Inc. (TSLA),  an American electric vehicle and clean energy company based in Palo Alto, California, for the two months of August and September of 2020 as listed in the $S\&P\, 500$ stock market index. We first estimate the model parameters using the method of moments.  

The expectation of $Y_t^j,$ return of SDE (\ref{equation1}), is given by 
\begin{eqnarray*} 
\mathbb{E}(Y_t^j)&=&r-\frac{1}{2} \sigma_j^2- \left(e^{\mu_{Z_{1,j,0}}+\frac{1}{2} \sigma_{Z_{1,j,0}}^2}-1\right) \lambda 
-\left(e^{\mu_{Z_{1,j,1}}+\frac{1}{2} \sigma_{Z_{1,j,1}}^2}-1\right) \lambda_j 
+\lambda\,\mu_{Z_{1,j,0}}\\ 
&&+\lambda_j\, \mu_{Z_{1,j,1}}. 
\end{eqnarray*} 
  
\begin{figure}[ht] 
\centerline{\includegraphics[height=6 cm, width=7 cm]{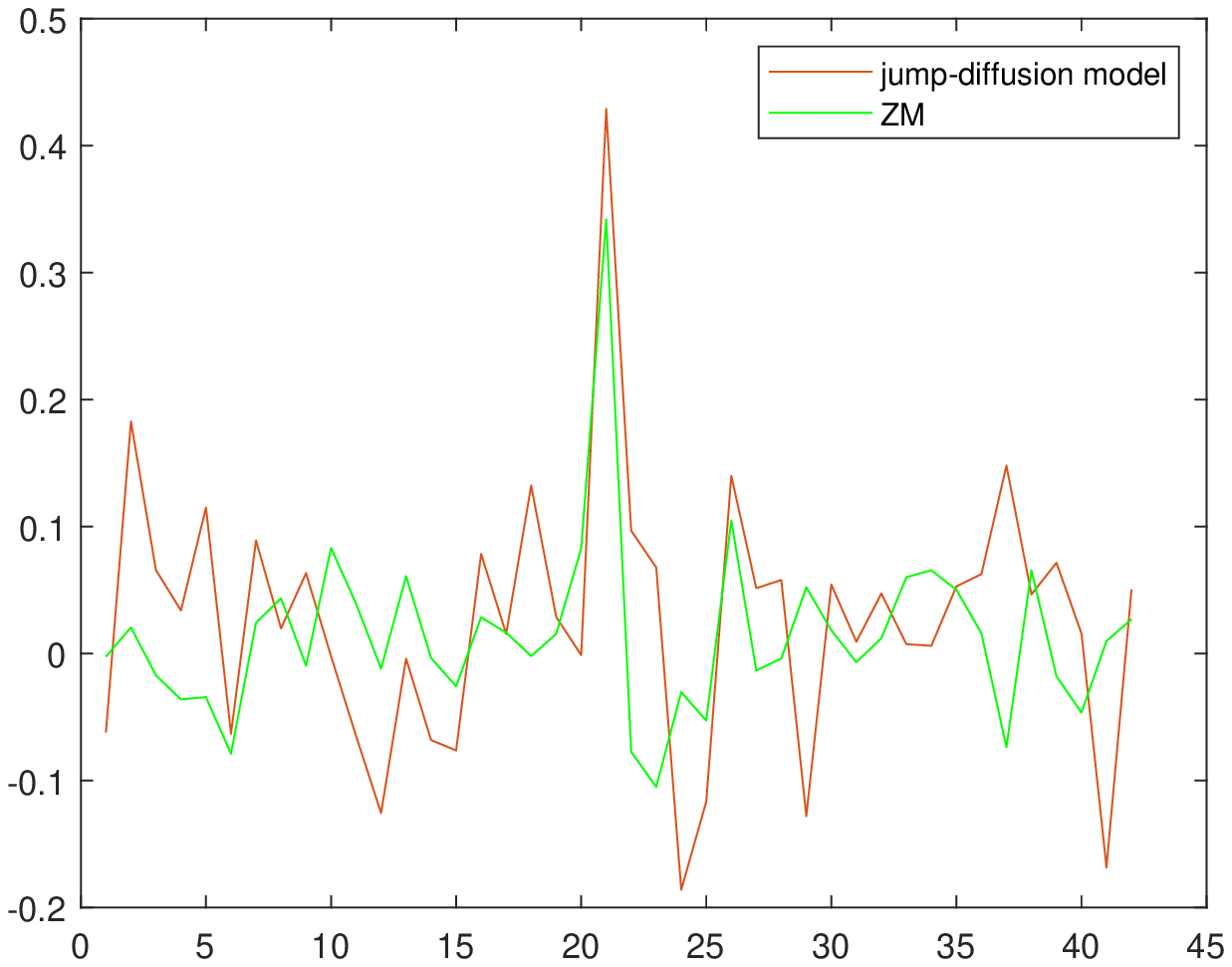} 
\includegraphics[height=6 cm, width=7 cm]{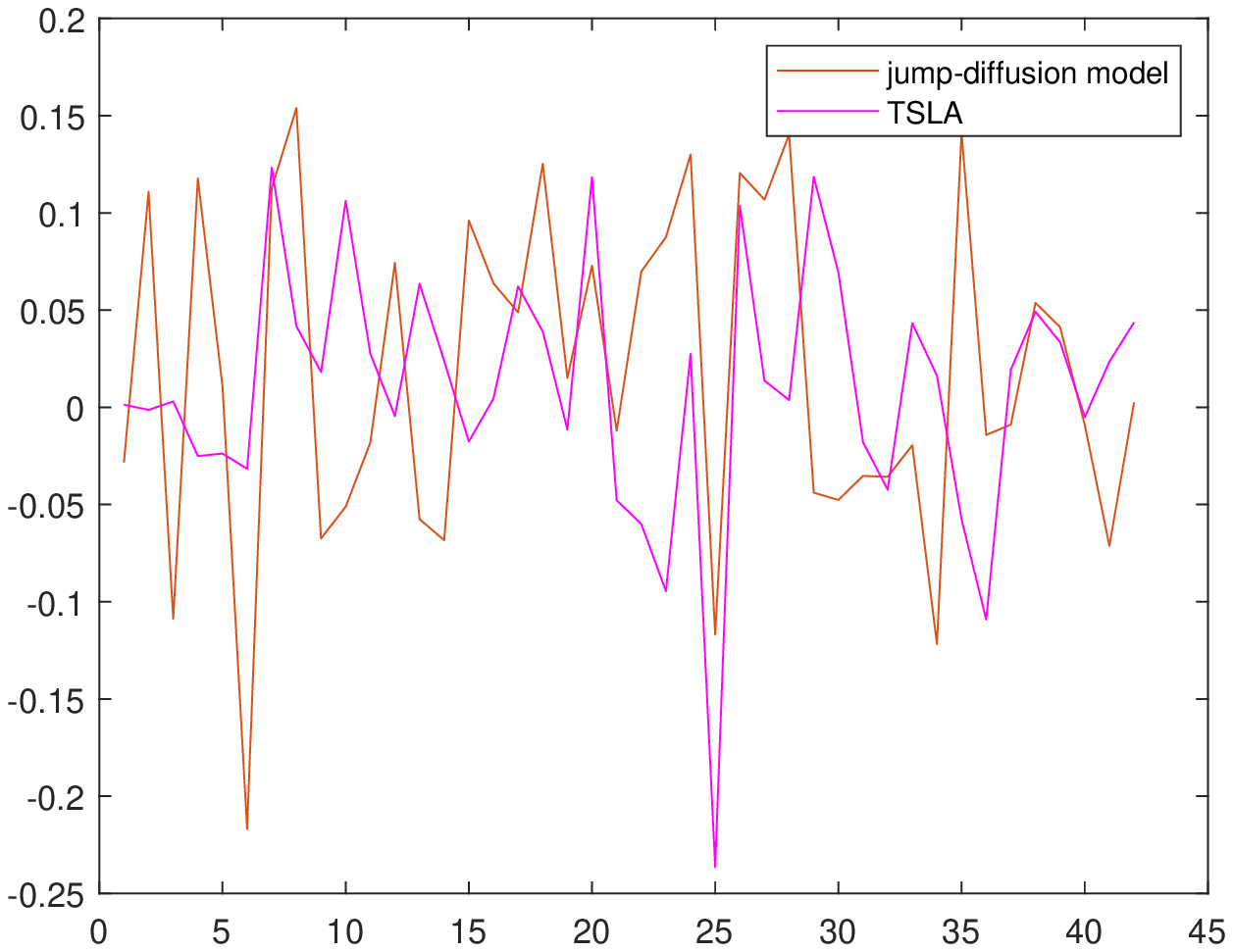}} 
\caption{\label{example5plot_ourmodel} Daily prices of ZM and TSLA stocks and their fitted values.} 
\end{figure}

We then use the Bootstrap Kolmogorov-Smirnov test to test whether the return of these assets follow Equation (\ref{equation1}). The $p$-values for ZM and TSLAS are respectively $0.112,\, 0.174,$ showing that the jump-diffusion model is a good fit for both data sets. Finally, we solve the optimization problem by the method described in Proposition \ref{final6}. The time horizon is $\tau=2$ months and we assume that the interest rate and the risk budget are $r=0.03$ and $p=0.05$ respectively. Figure \ref{example5plot_ourmodel} displays the actual prices together with their fitted values and Table \ref{mjy7688} illustrates that as in Example \ref{example3.11}, the terminal wealth decreases as $k^\star$ increases. 
 \begin{table} 
 \center \caption{Constant mix strategy.} \label{mjy7688} 
 \begin{tabular}{ c c c c c} 
    \noalign{\smallskip} \hline 
    $k^\star$ & $ (x_1,x_2)$   & $Terminal\, wealth$  & $Return$ \\  \hline  
    $0.5$ & $(0.4192,0.4922)$  & $2.5026$  & $22.3435 \%$ \\ 
    $0.85$ & $(0.3828,0.4495)$  & $2.4629$ & $20.4039 \%$  \\ 
     $0.9$ & $(0.2742,0.3219)$  & $2.3445$ & $14.6142 \%$  \\ 
    $0.95$ & $(0.1640,0.1926)$  & $1.3161$ & $8.7419 \%$  \\ \hline                               
 \end{tabular} 
 \end{table}

\eexam 
  
\begin{figure}[ht] 
\centerline{\includegraphics[height=6 cm, width=7 cm]{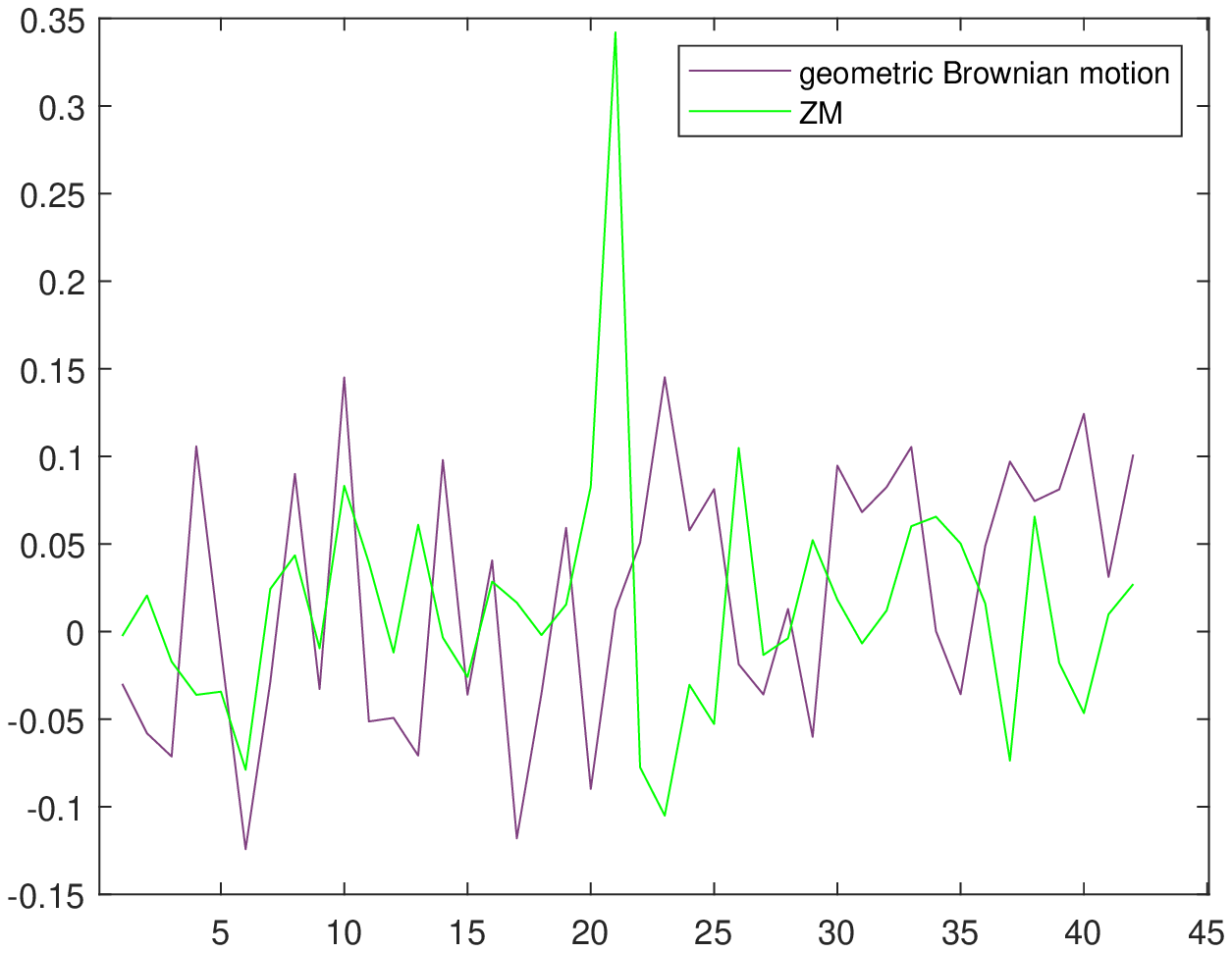} 
\includegraphics[height=6 cm, width=7 cm]{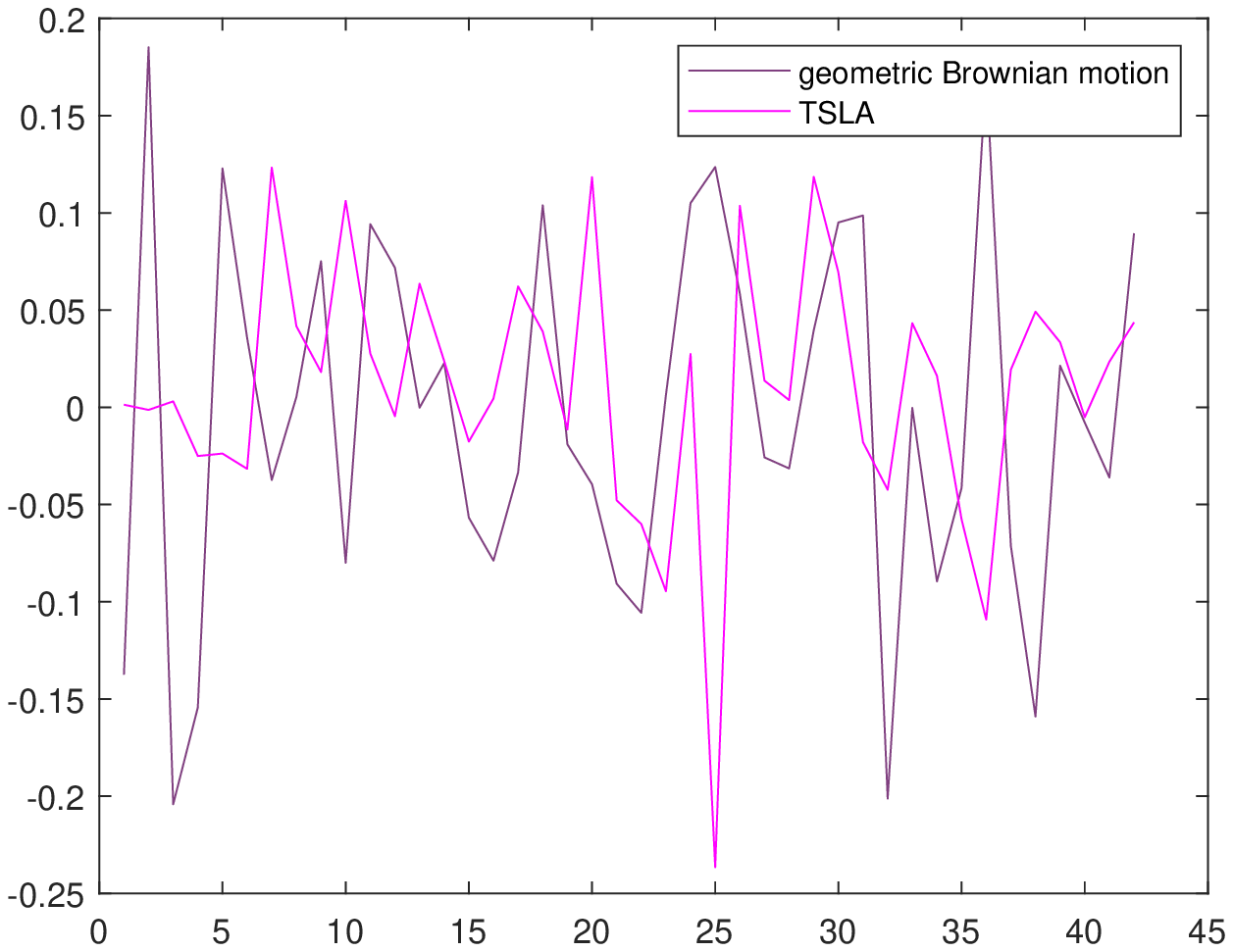}} 
\caption{\label{example5plot} Diagram of the daily prices of ZM and TSLA assets with the fitted values obtained form Geometric Brownian motion model.} 
\end{figure}

To assess the  performance of the proposed jump-diffusion model (\ref{equation1}) we shall compare it  with the geometric Brownian models. For this comparison we consider geometric Brownian motion in Dhaene et al. \cite{Dhaene:2005} with the following return equation, 
\begin{eqnarray}\label{jju55} 
\mathcal{Y}^j_t=(\kappa_j-\frac{1}{2} \gamma_j^2)+\gamma_j\, (B'_j(t)-B'_j(t-1)). 
\end{eqnarray} 
We first estimate the parameters of Equation (\ref{jju55}) using the method of moments and then run the corresponding Bootstrap Kolmogorov-Smirnov test (Figure \ref{example5plot}). The $p$-values for the return values of ZM and TSLA are respectively $0.09,\, 0.086,$ which are greater than the corresponding $p$-values obtained for our model which shows that the jump-diffusion model given in Equation (\ref{equation1}) has a better goodness of fit than the geometric Brownian model. 
\section{Conclusion and future work:} \label{Conclusion} 
In this paper, we studied portfolio optimization in a market with one risk-free and several risky assets for a constant mix investment strategy. We modeled the wealth process by a SDE with correlated jumps and designed a portfolio optimization problem based on controlling $CVaR.$ Since the solution to this optimization problem does not have a closed form, we obtained approximations using comonotonic bounds and solved the optimization problem for the terminal wealth. The performance of our method was assessed using a simulation study and an example from real data. Our ideas may be extended. For example, we may use a Levy process instead of the compound Poisson process, or a Levy copula can be implemented to model stock prices in a portfolio optimization problem. \\ 
  
{\bf{Acknowledgement:}} The authors are grateful to Professor Mahbanoo Tata for editing this paper as well as for her valuable suggestions. 
\bibliographystyle{amsplain} 
 

\newpage
\appendix 
\section{Appendix} 

{\bf{Proof of Proposition \ref{final1}: }} 
\begin{proof} 
To prove this proposition, it is enough to show that Equation (\ref{stochastic equation}) is the derivative of $P(t).$ We first  replace $P(t)$ with $f({\mathbf{G}}(t)).$ Let 
 ${\mathbf{G}}(t)=(G_1(t), G_2(t), G_3(t), G_4(t)),$ where 
$$G_1(t)=t,\,\,\,\, G_2(t)=B_{\mathbf{x}}(t),\,\,\,\, G_3(t)=\sum\limits_{j=1}^{m}\sum\limits_{k=1}^{N(t)}Z^\star_{k,j,0}({\mathbf{x}}),\,\,\,\, G_4(t)=\sum\limits_{j=1}^{m}\sum\limits_{k=1}^{N_j(t)}Z^\star_{k,j,1}({\mathbf{x}}),$$ 
and 
$$f({\mathbf{G}}(t))=e^{\left(\mu({\mathbf{x}})-\frac{{\mathbf{x}}^\top {\mathbf{\Sigma}}\, {\mathbf{x}}}{2} \right)\, G_1(t)+({\mathbf{x}}^\top {\mathbf{\Sigma}}\, {\mathbf{x}})\, G_2(t)+G_3(t)+G_4(t)}.$$ 
According to Ito's formula 
\begin{eqnarray*} 
f({\mathbf{G}}(t))-1&=&(\mu({\mathbf{x}})-\frac{{\mathbf{x}}^\top {\mathbf{\Sigma}}\, {\mathbf{x}}}{2} ) \int_{0^+}^{t}f({\mathbf{G}}(s-)) dG_1(s) 
+({\mathbf{x}}^\top {\mathbf{\Sigma}}\, {\mathbf{x}}) \int_{0^+}^{t}f({\mathbf{G}}(s-)) dG_2(s)\\ 
&&+\sum\limits_{i=3}^{4}\int_{0^+}^{t}f({\mathbf{G}}(s-))\, dG_i(s) 
+\frac{{\mathbf{x}}^\top {\mathbf{\Sigma}}\, {\mathbf{x}}}{2} \int_{0^+}^{t}f({\mathbf{G}}(s-))\, dG_1(s)\\ 
&&+\sum\limits_{0\leq s \leq t}\big\{[f({\mathbf{G}}(s))-f({\mathbf{G}}(s-))]-[f({\mathbf{G}}(s-))\, (\Delta G_3(s)+\Delta G_4(s))]\big\}. 
\end{eqnarray*} 
Simplifying, we have 
\begin{eqnarray*} 
f({\mathbf{G}}(t))-1&=&\mu({\mathbf{x}}) \int_{0^+}^{t}f({\mathbf{G}}(s-)) dG_1(s) 
+({\mathbf{x}}^\top {\mathbf{\Sigma}}\, {\mathbf{x}}) \int_{0^+}^{t}f({\mathbf{G}}(s-)) dG_2(s)\\ 
&&+\sum\limits_{i=3}^{4}\int_{0^+}^{t}f({\mathbf{G}}(s-))\, dG_i(s) 
+\sum\limits_{0\leq s \leq t} \left(f({\mathbf{G}}(s))-f({\mathbf{G}}(s-))\right)\\ 
&&-\sum\limits_{i=3}^{4} \sum\limits_{0\leq s \leq t}f({\mathbf{G}}(s-))\Delta G_i(s). 
\end{eqnarray*} 
For $i=3,4,$ we have  $\int_{0^+}^{t}f({\mathbf{G}}(s-)) dG_i(s)=\sum\limits_{0\leq s \leq t}f({\mathbf{G}}(s-))\, \Delta G_i(s),$ so that 
\begin{eqnarray*} 
f({\mathbf{G}}(t))-1&=&\mu({\mathbf{x}}) \int_{0^+}^{t}f({\mathbf{G}}(s-)) dG_1(s) 
+({\mathbf{x}}^\top {\mathbf{\Sigma}}\, {\mathbf{x}}) \int_{0^+}^{t}f({\mathbf{G}}(s-)) dG_2(s)\\ 
&&+\sum\limits_{0\leq s \leq t} \left(f({\mathbf{G}}(s))-f({\mathbf{G}}(s-))\right). 
\end{eqnarray*} 
Since 
\begin{eqnarray*} 
f({\mathbf{G}}(s))-f({\mathbf{G}}(s-))&=&e^{\left(\mu({\mathbf{x}})-\frac{{\mathbf{x}}^\top {\mathbf{\Sigma}}\, {\mathbf{x}}}{2} \right)\, G_1(s)+({\mathbf{x}}^\top {\mathbf{\Sigma}}\, {\mathbf{x}})\, G_2(s)+G_3(s)+G_4(s)}\\ 
&&-e^{\left(\mu({\mathbf{x}})-\frac{{\mathbf{x}}^\top {\mathbf{\Sigma}}\, {\mathbf{x}}}{2} \right)\, G_1(s-)+({\mathbf{x}}^\top {\mathbf{\Sigma}}\, {\mathbf{x}})\, G_2(s-)+G_3(s-)+G_4(s-)}\\ 
&=&e^{\left(\mu({\mathbf{x}})-\frac{{\mathbf{x}}^\top {\mathbf{\Sigma}}\, {\mathbf{x}}}{2} \right)\, G_1(s-)+({\mathbf{x}}^\top {\mathbf{\Sigma}}\, {\mathbf{x}})\, G_2(s-)+G_3(s-)+G_4(s-)}\\ 
&&\times e^{G_3(s)+G_4(s)-G_3(s-)-G_4(s-)}-1 \\ 
&=&f({\mathbf{G}}(s-))\left(e^{\Delta G_3(s)+\Delta G_4(s)}-1\right), 
\end{eqnarray*} 
hence 
\begin{eqnarray}\label{a11}\nonumber 
f({\mathbf{G}}(t))-1&=&\mu({\mathbf{x}}) \int_{0^+}^{t}f({\mathbf{G}}(s-))\, dG_1(s) 
+({\mathbf{x}}^\top {\mathbf{\Sigma}}\, {\mathbf{x}})\int_{0^+}^{t}f({\mathbf{G}}(s-))\, dG_2(s)\\ 
&&+\sum\limits_{0\leq s \leq t}f({\mathbf{G}}(s-))(e^{\Delta G_3(s)+\Delta G_4(s)}-1). 
\end{eqnarray} 
$G_{3}$ and $G_{4}$ do not have any common jumps, and so the summation in Equation (\ref{a11}) can be written as follows: 
\begin{eqnarray*} 
\sum\limits_{0\leq s \leq t}f({\mathbf{G}}(s-))(e^{\Delta G_3(s)}-1) 
+\sum\limits_{0\leq s \leq t}f({\mathbf{G}}(s-))(e^{\Delta G_4(s)}-1). 
\end{eqnarray*} 
Define
\begin{equation}\nonumber
Z^\star_{\iota(s-),j,0}({\mathbf{x}})=\left\{
\begin{array}{ll}
0, & \hbox{if\, there\, are\, no\, jumps\, in\, interval\, $(s-,s)$,} \\
Z^\star_{k,j,0}({\mathbf{x}}), & \hbox{if\, the\,$k$-th\, jump\,occurs \,in\, $(s-,s)$.}
\end{array}%
\right.
\end{equation} 
Now if the process $N$ does not have jumps at $s,$ then 
$\Delta G_{3}(s)=\sum\limits_{j=1}^{m} Z^\star_{\iota(s-),j,0}({\mathbf{x}})=0$ 
and 
$\Delta G_{4}(s)=\sum\limits_{j=1}^{m} Z^\star_{\iota(s-),j,1}({\mathbf{x}})=0,$ where

Therefore 
\begin{eqnarray*} 
&&\sum\limits_{0\leq s \leq t}f({\mathbf{G}}(s-))(e^{\Delta G_3(s)}-1)=\sum\limits_{0\leq s \leq t}f({\mathbf{G}}(s-))(e^{\sum\limits_{j=1}^{m} Z^\star_{\iota(s-),j,0}({\mathbf{x}})}-1)\, \Delta N(s) ,\\ 
&& \sum\limits_{0\leq s \leq t}f({\mathbf{G}}(s-))(e^{\Delta G_4(s)}-1)=\sum\limits_{0\leq s \leq t}f({\mathbf{G}}(s-))(e^{\sum\limits_{j=1}^{m} Z^\star_{\iota(s-),j,1}({\mathbf{x}})}-1)\, \Delta N_j(s), 
\end{eqnarray*}  
Then (\ref{a11}) reduces to 
\begin{eqnarray*} 
f({\mathbf{G}}(t))-1&=&\mu({\mathbf{x}}) \int_{0^+}^{t}f({\mathbf{G}}(s-))\, dG_1(s) 
+({\mathbf{x}}^\top {\mathbf{\Sigma}} {\mathbf{x}})\int_{0^+}^{t}f({\mathbf{G}}(s-))\, dG_2(s)\\ 
&&+\sum\limits_{0\leq s \leq t}f({\mathbf{G}}(s-))(e^{\sum\limits_{j=1}^{m} Z^\star_{\iota(s-),j,0}({\mathbf{x}})}-1)\, \Delta N(s) \\ 
&&+\sum\limits_{0\leq s \leq t}f({\mathbf{G}}(s-))(e^{\sum\limits_{j=1}^{m} Z^\star_{\iota(s-),j,1}({\mathbf{x}})}-1)\, \Delta N_j(s). 
\end{eqnarray*} 
For varying $j,\,$ $Z^\star_{j,0}({\mathbf{x}},s)$ do not have any common jumps and nor do $Z^\star_{j,1}({\mathbf{x}},s).$ Hence
\begin{eqnarray}\label{a12}  \nonumber 
f({\mathbf{G}}(t))-1&=&\mu({\mathbf{x}}) \int_{0^+}^{t}f({\mathbf{G}}(s-))\, dG_1(s) 
+({\mathbf{x}}^\top {\mathbf{\Sigma}} {\mathbf{x}})\int_{0^+}^{t}f({\mathbf{G}}(s-))\, dG_2(s)\\ \nonumber 
&&+\sum\limits_{j=1}^{m}\sum\limits_{0\leq s \leq t}f({\mathbf{G}}(s-))(e^{ Z^\star_{\iota(s-),j,0}({\mathbf{x}})}-1)\, \Delta N(s) \\ 
&&+\sum\limits_{j=1}^{m}\sum\limits_{0\leq s \leq t}f({\mathbf{G}}(s-))(e^{ Z^\star_{\iota(s-),j,1}({\mathbf{x}})}-1)\, \Delta N_j(s). 
\end{eqnarray} 
If $t_k$ is the time of the $k$-{th} jump, the first summation in (\ref{a12}) becomes 
$$\sum\limits_{j=1}^{m}\sum\limits_{0\leq s \leq t}f({\mathbf{G}}(s-))(e^{ Z^\star_{\iota(s-),j,0}({\mathbf{x}})}-1)=\sum\limits_{j=1}^{m}\sum\limits_{k=1}^{N(t)}f({\mathbf{G}}(t_k -))(e^{Z^\star_{k,j,0}({\mathbf{x}})}-1 ).$$ 
Now since $e^{Z^\star_{k,j,\ell}({\mathbf{x}})}-1=x_j(e^{Z_{k,j,\ell}}-1),\, \ell=0,1,$ so 
$$\sum\limits_{j=1}^{m}\sum\limits_{k=1}^{N(t)}f({\mathbf{G}}(t_k -))(e^{Z^\star_{k,j,0}({\mathbf{x}})}-1 )=\sum\limits_{j=1}^{m}\sum\limits_{k=1}^{N(t)}f({\mathbf{G}}(t_k -))x_j(e^{Z_{k,j,0}}-1 ).$$ 
Using a similar argument the second summation in (\ref{a12}) can be written as 
$$ \sum\limits_{j=1}^{m}\sum\limits_{k=1}^{N_j(t)}f({\mathbf{G}}(t_k -))x_j(e^{Z_{k,j,1}}-1 ). $$ 
 Thus (\ref{a12}) becomes 
\begin{eqnarray*} 
f({\mathbf{G}}(t))-1&=&\mu({\mathbf{x}}) \int_{0^+}^{t}f({\mathbf{G}}(s-))\, ds 
+({\mathbf{x}}^\top {\mathbf{\Sigma}}\, {\mathbf{x}})\int_{0^+}^{t}f({\mathbf{G}}(s-))\, dB_{\mathbf{x}}(s)\\ 
&&+\sum\limits_{j=1}^{m}\sum\limits_{k=1}^{N(t)}f({\mathbf{G}}(t_k -))x_j(e^{Z_{k,j,0}}-1 )\\ 
&&+\sum\limits_{j=1}^{m}\sum\limits_{k=1}^{N_j(t)}f({\mathbf{G}}(t_k -))x_j(e^{Z_{k,j,1}}-1 ). 
\end{eqnarray*} 
We have $\sum\limits_{k=1}^{N(t)}f({\mathbf{G}}(t_k -)) x_j (e^{Z_{k,j,0}}-1)=\int_{0^+}^{t}f({\mathbf{G}}(s -)) x_j (e^{Z_{\iota(s-),j,0}}-1 )\, dN(s)$ and 
  
 $\sum\limits_{k=1}^{N_j(t)}f({\mathbf{G}}(t_k -)) x_j (e^{Z_{k,j,1}}-1 )=\int_{0^+}^{t} f({\mathbf{G}}(s -)) x_j (e^{Z_{\iota(s-),j,1}}-1 )\, dN_j(s),$ so 
\begin{eqnarray*} 
f({\mathbf{G}}(t))-1&=&\mu({\mathbf{x}}) \int_{0^+}^{t}f({\mathbf{G}}(s-))\, ds 
+({\mathbf{x}}^\top {\mathbf{\Sigma}}\, {\mathbf{x}})\int_{0^+}^{t}f({\mathbf{G}}(s-))\, dB_{\mathbf{x}}(s)\\ 
&&+\int_{0^+}^{t}f({\mathbf{G}}(s -)) x_j (e^{Z_{\iota(s-),j,0}}-1 )\, dN(s)\\ 
&&+\int_{0^+}^{t} f({\mathbf{G}}(s -)) x_j (e^{Z_{\iota(s-)j,1}}-1 )\, dN_j(s). 
\end{eqnarray*} 
Differentiating, 
\begin{eqnarray*} 
\frac{d f({\mathbf{G}}(t))}{d t}&=&\mu({\mathbf{x}})f({\mathbf{G}}(t-))\, dt 
+({\mathbf{x}}^\top {\mathbf{\Sigma}}\, {\mathbf{x}}) f({\mathbf{G}}(t-))\, dB_{\mathbf{x}}(t) 
+f({\mathbf{G}}(t -))\\ 
&&\times \sum\limits_{j=1}^{m}x_j\, (e^{Z_{\iota(t-),j,0}}-1 )\, dN(t) 
+ f({\mathbf{G}}(t-))\sum\limits_{j=1}^{m}x_j\,(e^{Z_{\iota(t-),j,1}}-1 )\, dN_j(t). 
\end{eqnarray*} 
  
\begin{eqnarray*} 
\frac{d f({\mathbf{G}}(t))}{d t}&=&\mu({\mathbf{x}})f({\mathbf{G}}(t-))\, dt 
+({\mathbf{x}}^\top {\mathbf{\Sigma}}\, {\mathbf{x}}) f({\mathbf{G}}(t-))\, dB_{\mathbf{x}}(t) 
+f({\mathbf{G}}(t -))\\ 
&&\times \sum\limits_{j=1}^{m}x_j\, (e^{Z_{\iota(t-),j,0}}-1 )\, dN(t) 
+ f({\mathbf{G}}(t-))\sum\limits_{j=1}^{m}x_j\,(e^{Z_{\iota(t-),j,1}}-1 )\, dN_j(t). 
\end{eqnarray*} 
Dividing both sides of the above equation by $f({\mathbf{G}}(t-)),$ we obtain (\ref{stochastic equation}). By Theorem 6  [\cite{Protter}, P 249 ], the solution is unique. 
\end{proof} 

{\bf{Proof of Theorem \ref{Gnllllnnn}: }} 
\begin{proof} 
In what follows, $M_Z(\theta)=\mathbb{E}(e^{Z\theta})$ presents the moment generating function of the random variable $Z\,.$
We calculate the lower bound for terminal wealth using Theorem \ref{theorem1.4}. 
\begin{eqnarray*} 
W_\tau^L&=&\sum\limits_{t=0}^{\tau-1}\alpha_t\, \mathbb{E}(e^{S_t})\, \mathbb{E}(e^{V_t}\mid \Lambda) 
=\sum\limits_{t=0}^{\tau-1}\alpha_t\, e^{(\tau-t)(\mu({\mathbf{x}})-\frac{{\mathbf{x}}^\top {\mathbf{\Sigma}}\, {\mathbf{x}}}{2})} \\ 
&&\times \mathbb{E}\left(e^{\sum\limits_{j=1}^{m}\sum\limits_{k=N(t)+1}^{N(\tau)}Z^\star_{k,j,0}({\mathbf{x}})}\right) \,\mathbb{E}\left(e^{\sum\limits_{j=1}^{m}\sum\limits_{k=N_j(t)+1}^{N_j(\tau)}Z^\star_{k,j,1}({\mathbf{x}})}\right) 
\, M_{V_t\mid \Lambda}(1). 
\end{eqnarray*} 
We know that 
$V_t\mid \Lambda =\lambda \sim N\left(r_t \, \sigma_{V_t}\, \frac{\lambda}{\sigma_{\Lambda}}, \sigma_{V_t}^2\, (1-r_t^2) \right), $ 
where $r_t=Corr(V_t,\Lambda).$ Hence 
$M_{V_t\mid \Lambda}(1)=e^{r_t\, \sigma_{V_t}\, \frac{\Lambda}{\sigma_\Lambda}+\frac{1}{2}\sigma_{V_t}^2\, (1-r_t^2) }.$ 
  
On the other hand, 
\begin{eqnarray*} 
&&\mathbb{E}\left(e^{\sum\limits_{j=1}^{m}\sum\limits_{k=N(t)+1}^{N(\tau)}Z^\star_{k,j,0}({\mathbf{x}})}\right)=e^{\lambda\,(\tau-t)\,\left(\prod\limits_{j=1}^{m} M_{Z^\star_{1,j,0}({\mathbf{x}})}(1)-1\right)}. 
\end{eqnarray*} 
Similarly, 
\begin{eqnarray*} 
&&\mathbb{E}\left(e^{\sum\limits_{j=1}^{m}\sum\limits_{k=N_j(t)+1}^{N_j(\tau)}Z^\star_{k,j,1}({\mathbf{x}})}\right) 
= e^{-(\tau-t)\sum\limits_{j=1}^{m} \lambda_j \left(M_{Z^\star_{j,1}({\mathbf{x}})}(1)-1\right)}. 
\end{eqnarray*} 
Therefore 
\begin{eqnarray*} 
W_\tau^L&=&\sum\limits_{t=0}^{\tau-1}\alpha_t\,\, e^{(\tau-t)(\mu({\mathbf{x}})-\frac{{\mathbf{x}}^\top {\mathbf{\Sigma}}\, {\mathbf{x}}}{2})}\,\, e^{r_t\, \sigma_{V_t}\, \frac{\Lambda}{\sigma_\Lambda}+\frac{1}{2}\sigma_{V_t}^2\, (1-r_t^2) }\,\,\\ 
&&\times e^{\lambda\,(\tau-t)\left(\prod_{j=1}^{m}M_{Z^\star_{1,j,0}({\mathbf{x}})}(1)-1\right) 
+\sum_{j=1}^{m} \lambda_j\,(\tau-t)\left(M_{Z^\star_{1,j,1}({\mathbf{x}})}(1)-1\right)}. 
\end{eqnarray*} 
We now calculate $\sigma^2_{V_t}$ and $r_t.$ By definition, $B_{{\mathbf{x}}}(t):=\frac{1}{{\mathbf{x}}^\top\,{\mathbf{\Sigma}}\,\, {\mathbf{x}}}\sum\limits_{j=1}^{m} x_j\, \sigma_j\, B'_j(t),$ in which $B'_i$ are standard Brownian motions with $Cov(B'_i(t),B'_j(t+s))=t\,\rho_{i,j}\, (t,s\geq 0,\,   i,j=1,\ldots,m)\,,$ therefore $B_{{\mathbf{x}}}(t)\sim N(0,\frac{t}{{\mathbf{x}}^\top {\mathbf{\Sigma}}\, {\mathbf{x}}}) $ and
\begin{eqnarray*} 
\sigma^2_{V_t}&=&({\mathbf{x}}^\top {\mathbf{\Sigma}}\, {\mathbf{x}})^2 \,Var\left(B_{{\mathbf{x}}}(\tau-t) \right) 
=(\tau-t)\, ({\mathbf{x}}^\top {\mathbf{\Sigma}}\, {\mathbf{x}}). 
\end{eqnarray*} 
Also $r_t =Corr(V_t,\Lambda)=\frac{Cov(V_t,\Lambda)}{\sigma_{V_t}\, \sigma_{\Lambda}},$ so we need to obtain $\sigma_{\Lambda}$ and $Cov(V_t,\Lambda).$ First, \newline 
$\sigma^2_{\Lambda}=Var\left(\sum\limits_{t=0}^{\tau-1}\alpha_t\, V_t \right) 
=Cov\left(\sum\limits_{t=0}^{\tau-1}\alpha_t\, V_t, \sum\limits_{l=0}^{\tau-1}\alpha_l\, V_l  \right),$ so 
\begin{eqnarray*} 
\sigma^2_{\Lambda}&=&({\mathbf{x}}^\top {\mathbf{\Sigma}}\, {\mathbf{x}})^2\, \sum\limits_{t,l=0}^{\tau-1} \alpha_t\,\alpha_l 
\big\{Cov(B_{{\mathbf{x}}}(\tau),B_{{\mathbf{x}}}(\tau)) 
-Cov(B_{{\mathbf{x}}}(\tau),B_{{\mathbf{x}}}(l)) \\ 
&&-Cov(B_{{\mathbf{x}}}(t),B_{{\mathbf{x}}}(\tau)) 
+Cov(B_{{\mathbf{x}}}(t),B_{{\mathbf{x}}}(l)) \big\}\\ 
&=&({\mathbf{x}}^\top {\mathbf{\Sigma}}\, {\mathbf{x}})^2\, \sum\limits_{t,l=0}^{\tau-1}\alpha_t\,\alpha_l 
\big\{\frac{\tau}{{\mathbf{x}}^\top {\mathbf{\Sigma}}\, {\mathbf{x}}}-\frac{l}{{\mathbf{x}}^\top {\mathbf{\Sigma}}\, {\mathbf{x}}}-\frac{t}{{\mathbf{x}}^\top {\mathbf{\Sigma}}\, {\mathbf{x}}}+\frac{\min\{t,l\}}{{\mathbf{x}}^\top {\mathbf{\Sigma}}\, {\mathbf{x}}} \big\}\\ 
&=&({\mathbf{x}}^\top {\mathbf{\Sigma}}\, {\mathbf{x}})\, \sum\limits_{t,l=0}^{\tau-1}\alpha_t\,\alpha_l \min\{\tau-t,\tau-l\}. 
\end{eqnarray*} 
Similarly $Cov(V_t,\Lambda)= ({\mathbf{x}}^\top {\mathbf{\Sigma}}\, {\mathbf{x}}) \sum\limits_{l=0}^{\tau-1}\alpha_l\,\min\{\tau-t,\tau-l\},$ 
hence $r_t=\frac{ \sum\limits_{l=0}^{\tau-1}\alpha_l\,\min\{\tau-t,\tau-l\} } 
{\sqrt{(\tau-t)  \sum\limits_{k,l=0}^{\tau-1}\alpha_k\, \alpha_l\,\min\{\tau-k,\tau-l\}}}$ and 
\begin{eqnarray*} 
W_\tau^L&=&\sum\limits_{t=0}^{\tau-1}\alpha_t\, e^{(\tau-t) ({\mathbf{\mu}}-r {\mathbf{1}})^\top {\mathbf{x}}+(\tau-t)r-\frac{1}{2}(\tau-t)\, ({\mathbf{x}}^\top {\mathbf{\Sigma}}\, {\mathbf{x}})+\frac{ \sum\limits_{l=0}^{\tau-1}\alpha_l\,\min\{\tau-t,\tau-l\} } 
{\sqrt{ \sum\limits_{k,l=0}^{\tau-1}\alpha_k\, \alpha_l\,\min\{\tau-k,\tau-l\} }}\, \sqrt{{\mathbf{x}}^\top {\mathbf{\Sigma}}\, {\mathbf{x}}}\,\frac{\Lambda}{\sigma_\Lambda} }\\ 
&&\times e^{\frac{1}{2}(\tau-t)\, ({\mathbf{x}}^\top {\mathbf{\Sigma}}\, {\mathbf{x}})-\frac{1}{2}\,({\mathbf{x}}^\top {\mathbf{\Sigma}}\, {\mathbf{x}}) 
\frac{ \left(\sum\limits_{l=0}^{\tau-1}\alpha_l\,\min\{\tau-t,\tau-l\}\right)^2 } 
{\sum\limits_{k,l=0}^{\tau-1}\alpha_k\, \alpha_l\,\min\{\tau-k,\tau-l\} }}\\ 
&&\times e^{\lambda\,(\tau-t)\left(\prod_{j=1}^{m}M_{Z^\star_{1,j,0}({\mathbf{x}})}(1)-1\right) 
+\sum_{j=1}^{m} \lambda_j\,(\tau-t)\left(M_{Z^\star_{1,j,1}({\mathbf{x}})}(1)-1\right)}. 
\end{eqnarray*} 
We may write 
\begin{eqnarray*} 
W_\tau^L=\sum\limits_{t=0}^{\tau-1}\alpha_t\, e^{\frac{\Lambda}{\sigma_\Lambda}\, \sqrt{{\mathbf{x}}^\top {\mathbf{\Sigma}}\, {\mathbf{x}}}\,\, c_{4,t}+c_{3,t}({\mathbf{x}})+c_{2,t}}, 
\end{eqnarray*} 
where 
\begin{eqnarray}\label{a14} \nonumber
 c_{2,t}&=&(\tau-t)r+\lambda\,(\tau-t)\left(\prod_{j=1}^{m}M_{Z^\star_{1,j,0}({\mathbf{x}})}(1)-1\right) \\
&&+\sum_{j=1}^{m} \lambda_j\,(\tau-t)\left(M_{Z^\star_{1,j,1}({\mathbf{x}})}(1)-1\right), 
\end{eqnarray} 
\begin{eqnarray}\label{a15} 
c_{3,t}({\mathbf{x}})=(\tau-t) ({\mathbf{\mu}}-r {\mathbf{1}})^\top {\mathbf{x}}-\frac{1}{2}\, ({\mathbf{x}}^\top {\mathbf{\Sigma}}\, {\mathbf{x}})\, 
\frac{\left (\sum\limits_{l=0}^{\tau-1}\alpha_l\,\min\{\tau-t,\tau-l\}\right)^2 } 
{\sum\limits_{k,l=0}^{\tau-1}\alpha_k\, \alpha_l\,\min\{\tau-k,\tau-l\}}, 
\end{eqnarray} 
and 
\begin{eqnarray}\label{a16} 
 \hspace{-5.2cm} c_{4,t}=\frac{ \sum\limits_{l=0}^{\tau-1}\alpha_l\,\min\{\tau-t,\tau-l\} } 
{\sqrt{\sum\limits_{k,l=0}^{\tau-1}\alpha_k\, \alpha_l\,\min\{\tau-k,\tau-l\}}}, 
\end{eqnarray} 
which completes the result. 
\end{proof} 

{\bf{Proof of Corollary \ref{corollary3.77}: }} 
\begin{proof} 
From (\ref{termG_n^lnnnnn}) 
\begin{eqnarray*} 
W_\tau^{'L}&=&\sum\limits_{t=0}^{\tau-1} \alpha_t\, \left(1+\frac{\Lambda}{\sigma_\Lambda}\, \sqrt{{\mathbf{x}}^\top {\mathbf{\Sigma}}\, {\mathbf{x}}}\,\, c_{4,t}+c_{3,t}({\mathbf{x}})+c_{2,t} \right)\\ 
&=&\sum\limits_{t=0}^{\tau-1} \alpha_t\, (1+c_{2,t})+\sum\limits_{t=0}^{\tau-1} \alpha_t\, c_{3,t}({\mathbf{x}}) 
+\frac{\Lambda}{\sigma_\Lambda}\, \sqrt{{\mathbf{x}}^\top {\mathbf{\Sigma}}\, {\mathbf{x}}}\,\sum\limits_{t=0}^{\tau-1} \alpha_t\, c_{4,t}\\ 
&=&\sum\limits_{t=0}^{\tau-1} \alpha_t\, (1+c_{2,t}) 
+({\mathbf{\mu}}-r {\mathbf{1}})^\top {\mathbf{x}}\,\, \left(\sum\limits_{t=0}^{\tau-1} \alpha_t\,(\tau-t) \right)\\ 
&&-({\mathbf{x}}^\top {\mathbf{\Sigma}}\, {\mathbf{x}})\,\sum\limits_{t=0}^{\tau-1} \alpha_t\, \frac{\left (\sum\limits_{l=0}^{\tau-1}\alpha_l\, 
\min\{\tau-t,\tau-l\}\right)^2 } {2\sum\limits_{k,l=0}^{\tau-1}\alpha_k\, \alpha_l\,\min\{\tau-k,\tau-l\}}\\ 
&&+\frac{\Lambda}{\sigma_\Lambda}\, \sqrt{{\mathbf{x}}^\top {\mathbf{\Sigma}}\, {\mathbf{x}}}\,\sum\limits_{t=0}^{\tau-1} \alpha_t\, 
\frac{ \sum\limits_{l=0}^{\tau-1}\alpha_l\,\min\{\tau-t,\tau-l\} } 
{\sqrt{\sum\limits_{k,l=0}^{\tau-1}\alpha_k\, \alpha_l\,\min\{\tau-k,\tau-l\}}}. 
\end{eqnarray*} 
So 
\begin{eqnarray*} 
W_\tau^{'L}=c_5+c_6\,\, ({\mathbf{\mu}}-r {\mathbf{1}})^\top {\mathbf{x}} -c_7\,({\mathbf{x}}^\top {\mathbf{\Sigma}}\, {\mathbf{x}})+\frac{\Lambda}{\sigma_\Lambda}\,c_8\, \sqrt{{\mathbf{x}}^\top {\mathbf{\Sigma}}\, {\mathbf{x}}}, 
\end{eqnarray*} 
where 
\begin{eqnarray}\label{a17} 
\hspace{-5.1cm} c_5=\sum\limits_{t=0}^{\tau-1} \alpha_t\, (1+c_{2,t}), 
\end{eqnarray} 
\begin{eqnarray}\label{a18} 
\hspace{-5.5cm} c_6=\sum\limits_{t=0}^{\tau-1} \alpha_t\,(\tau-t), 
\end{eqnarray} 
\begin{eqnarray}\label{a19} 
\hspace{-1.8cm} c_7=\sum\limits_{t=0}^{\tau-1} \alpha_t\, \frac{\left (\sum\limits_{l=0}^{\tau-1}\alpha_l\,\min\{\tau-t,\tau-l\}\right)^2 } {2\sum\limits_{k,l=0}^{\tau-1}\alpha_k\, \alpha_l\,\min\{\tau-k,\tau-l\}}, 
\end{eqnarray} 
\begin{eqnarray}\label{a20} 
\hspace{-1.6cm} c_8=\sum\limits_{t=0}^{\tau-1} \alpha_t\, 
\frac{ \sum\limits_{l=0}^{\tau-1}\alpha_l\,\min\{\tau-t,\tau-l\} } 
{\sqrt{\sum\limits_{k,l=0}^{\tau-1}\alpha_k\, \alpha_l\,\min\{\tau-k,\tau-l\}}}, 
\end{eqnarray} 
this shows the result. 
\end{proof} 
  \newpage
{\bf{Proof of Lemma \ref{gghhjk54}: }} 
\begin{proof} 
From (\ref{h34567}) 
\begin{eqnarray*} 
CVaR_{1-p}(-W_\tau^{'L})&=&CVaR_{1-p}\left(-c_{5}- c_{6}\, (({\mathbf{\mu}}-r {\mathbf{1}})^\top\, {\mathbf{x}})+c_7\, ({\mathbf{x}}^\top {\mathbf{\Sigma}}\, {\mathbf{x}})-c_8\, \frac{\Lambda}{\sigma_\Lambda}\, \sqrt{{\mathbf{x}}^\top {\mathbf{\Sigma}}\, {\mathbf{x}}}\right)\\ 
&=&-c_{5}- c_{6}\,\, ({\mathbf{\mu}}-r {\mathbf{1}})^\top\, {\mathbf{x}}+c_7\, ({\mathbf{x}}^\top {\mathbf{\Sigma}}\, {\mathbf{x}})+c_8\,\sqrt{{\mathbf{x}}^\top {\mathbf{\Sigma}}\, {\mathbf{x}}}\,\, CVaR_{1-p} (- \frac{\Lambda}{\sigma_\Lambda}). 
\end{eqnarray*} 
Let $Z=\frac{\Lambda}{\sigma_\Lambda}$ 
\begin{eqnarray*} 
CVaR_{1-p} (-Z)&=&\mathbb{E} [-Z\mid  -Z>VaR_{1-p}(-Z)] 
=\mathbb{E} [-Z\mid  -Z>-\Phi^{-1}(p) ]\\ 
&=&-\mathbb{E} [Z \mid  Z<\Phi^{-1}(p) ] 
=-\int\limits_{-\infty}^{+\infty}z\, f_{Z \mid  Z<\Phi^{-1}(p) }(z)\, dz\\ 
&=&-\int\limits_{-\infty}^{+\infty}z\, \frac{f_{Z ,  Z<\Phi^{-1}(p) }(z)}{P(Z<\Phi^{-1}(p))} 
=-\frac{1}{p}\int\limits_{-\infty}^{\Phi^{-1}(p)}z\, f_Z(z)\,dz,\\ 
\end{eqnarray*} 
where the r.v., $Z$ has standard normal distribution, 
\begin{eqnarray*} 
CVaR_{1-p} (-Z)=-\frac{1}{p}\int\limits_{-\infty}^{\Phi^{-1}(p)}z\,\, \frac{1}{\sqrt{2\, \pi}}\,e^{-\frac{1}{2}z^2} dz 
=\frac{1}{p\, \sqrt{2\, \pi}} e^{-\frac{1}{2}(\Phi^{-1}(p))^2}. 
\end{eqnarray*} 
Thus 
\begin{eqnarray*} 
\hspace{.3cm} CVaR_{1-p}(-W_\tau^{'L})=-c_{5}- c_{6}\,\, ({\mathbf{\mu}}-r {\mathbf{1}})^\top\, {\mathbf{x}}+c_7\, ({\mathbf{x}}^\top {\mathbf{\Sigma}}\, {\mathbf{x}})+c_9\, \sqrt{{\mathbf{x}}^\top {\mathbf{\Sigma}}\, {\mathbf{x}}}, 
\end{eqnarray*} 
where 
\begin{eqnarray}\label{a21} 
\hspace{-6cm} c_9=c_8\,\frac{1}{p\, \sqrt{2\, \pi}} e^{-\frac{1}{2}(\Phi^{-1}(p))^2}, 
\end{eqnarray} 
and this completes the proof of the result. 
\end{proof} 

{\bf{Proof of Proposition \ref{final4}: }} 
\begin{proof} 
We optimize the portfolio using the fractional Kelly-strategy. For this purpose, consider Lagrange's equation 
\begin{eqnarray*} 
L({\mathbf{x}}, \lambda_1, \lambda_2)&=& c_5+c_6\,\, ({\mathbf{\mu}}-r {\mathbf{1}})^\top {\mathbf{x}}-c_7 ({\mathbf{x}}^\top {\mathbf{\Sigma}}\, {\mathbf{x}}) 
+\lambda_1 \big(-c_{5}- c_{6}\,\, ({\mathbf{\mu}}-r {\mathbf{1}})^\top {\mathbf{x}}+c_7 ({\mathbf{x}}^\top {\mathbf{\Sigma}}\, {\mathbf{x}})\\ 
&&+c_9 \sqrt{{\mathbf{x}}^\top {\mathbf{\Sigma}}\, {\mathbf{x}}}+K \big) 
+\lambda_2 \big(r+ ({\mathbf{\mu}}-r {\mathbf{1}})^\top {\mathbf{x}}-c_0 \big). 
\end{eqnarray*} 
where $ {\mathbf{\mu}}:=(r+\mu_1-\lambda h_{1,0}-\lambda_1 h_{1,1},\ldots,r+\mu_m-\lambda h_{m,0}-\lambda_m h_{m,1})^\top,$ and $K$ is as in (\ref{n88}). 
It can be easily shown that the derivative of $L({\mathbf{x}}, \lambda_1, \lambda_2)$ is as follows: 
\begin{eqnarray*} 
c_6\, ({\mathbf{\mu}}-r {\mathbf{1}})-2\, c_7\, ({\mathbf{\Sigma}}\, {\mathbf{x}})-\lambda_1\, c_{6}\,\, ({\mathbf{\mu}}-r {\mathbf{1}})+2\, \lambda_1\,c_7\, ({\mathbf{\Sigma}}\, {\mathbf{x}})+\lambda_1\,c_9 
\frac{{\mathbf{\Sigma}}\, {\mathbf{x}}}{\sqrt{{\mathbf{x}}^\top {\mathbf{\Sigma}}\, {\mathbf{x}}}}+\lambda_2({\mathbf{\mu}}-r {\mathbf{1}})={\mathbf{0}}. 
\end{eqnarray*} 
This implies that 
\begin{eqnarray*} 
{\mathbf{x}}={{\mathbf{\Sigma}}^{-1}\, ({\mathbf{\mu}}-r {\mathbf{1}})}\,\, \frac{c_6\, (\lambda_1-1)-\lambda_2}{2\, c_7\,(\lambda_1-1)+\frac{\lambda_1\,c_9 }{\sqrt{{\mathbf{x}}^\top {\mathbf{\Sigma}}\, {\mathbf{x}}}}}. 
\end{eqnarray*} 
Thus ${\mathbf{x}}={\mathbf{x}}^\star\, q.$ 
\end{proof} 


  
\end{document}